\documentclass[notitlepage,leqno,10pt]{article}
\usepackage{amsmath, amsfonts, amssymb, amscd,
amsthm,enumitem,extarrows}
\usepackage[all]{xy}
\usepackage{latexsym}
\usepackage{amsmath}
\usepackage{amsfonts}
\usepackage{amssymb}
\usepackage{amscd}
\usepackage{mathtools, amsthm,enumitem}

\newcommand{\fa}{\forall}

\renewcommand{\d}{\delta } 
\newcommand{\D }{\Delta }

\renewcommand{\l }{\lambda }

\newcommand{\ov}{\overline}
\newcommand{\intbar}{\mathop{\int\makebox(-13.5,0){\rule[4pt]{.7em}{0.3pt}}%
\kern-6pt}\nolimits}

\newcommand{\be}{\begin{equation}}
\newcommand{\ee}{\end{equation}}
\newcommand{\bes}{\begin{equation*}}
\newcommand{\ees}{\end{equation*}}
\newcommand{\ba}{\begin{eqnarray}}
\newcommand{\ea}{\end{eqnarray}}
\newcommand{\bas}{\begin{eqnarray*}}
\newcommand{\eas}{\end{eqnarray*}}
\newenvironment{pf}{\noindent{\sc Proof}.\enspace}{\rule{2mm}{2mm}\medskip}
\newenvironment{pfn}{\noindent{\sc Proof}}{\rule{2mm}{2mm}\medskip}

\newcommand{\R}{\mathbb{R}}

\newcommand{\N}{\mathbb{N}}

\author{Martin MAYER$^{a}$,\;\; Cheikh Birahim NDIAYE$^{b}$}
\date{}

\title{\bf Asymptotics of the Poisson kernel and Green's functions of the fractional conformal Laplacian}

\begin{document}

\newtheorem{lem}{Lemma}[section]
\newtheorem{pro}[lem]{Proposition}
\newtheorem{thm}[lem]{Theorem}
\newtheorem{rem}[lem]{Remark}
\newtheorem{cor}[lem]{Corollary}
\newtheorem{df}[lem]{Definition}

\maketitle

\begin{center}

{\small

\noindent  $^{a}$ Dipartimento di Matematica della Universit\'a di Roma Tor Vergata, \\ 
Via della Ricerca Scientifica 1, 00133 Roma, ITALY.

}
\
\
\
\
 
{\small

\noindent  $^{b}$ Department of Mathematics  of Howard University \\  Annex 3, Graduate School of Arts and Sciences, \# 217 \\ DC 20059 Washington, USA.

}

\
\
{\small

\noindent

}

\end{center}

\footnotetext[1]
{
E-mail: mayer@mat.uniroma2.it,  
cheikh.ndiaye@howard.edu.
\\
\thanks
{
\\ C. B. Ndiaye was partially supported by NSF grant DMS--2000164
\\  M.Mayer has been supported by the Italian MIUR Department of Excellence grant CUP E83C18000100006. 
}
}

\

\

\begin{center}
{\bf Abstract}

\end{center}
We study the asymptotics of the Poisson kernel and Green's functions of the fractional conformal Laplacian for conformal infinities of asymptotically hyperbolic manifolds. We derive sharp expansions of the Poisson kernel and Green's functions of the conformal Laplacian near their singularities. Our expansions of the Green's functions answer the first part of the conjecture of Kim-Musso-Wei\cite{kmw1} in the case of locally flat conformal infinities of Poincare-Einstein manifolds and together with the Poisson kernel asymptotic is used also in our paper \cite{martndia2} to show solvability of the fractional Yamabe problem in that case. Our asymptotics of the Green's functions on the general case of conformal infinities of asymptotically hyperbolic space is used also in \cite{nss} to show solvability of the fractional Yamabe problem for conformal infinities of dimension \;$3$ and fractional parameter in \;$(\frac{1}{2}, 1)$\; to a global case left by previous works.

\begin{center}
 \
 \


\bigskip

\centerline{\bf AMS subject classification:  53C21, 35C60, 58J60, 55N10.}

\bigskip
\end{center}

\begin{center}
$_{}$
\vspace{100pt}
\tableofcontents
\end{center}
\section{Introduction and statement of the results}
In the last decades there has been a lot of study about fractional order operators in Analysis and Geometric Analysis as well.  In both fields, the recurrent themes are existence, regularity and sharp estimates,  see \cite{cafroq}, \cite{cafsyl}, \cite{cafsoug}, \cite{cafval}, \cite{cafvas}, \cite{cg}, \cite{fg},  \cite{fabkenser}, \cite{qr}, \cite{gz}, \cite{guil}, \cite{gms}, \cite{gq}). In this paper we are interested in the issue of existence, regularity and sharp estimates in the context of Conformal Geometry. Precisely, we study the issue of existence, regularity and sharp asymptotics of the Poisson and Green's functions of the fractional conformal Laplacian on conformal infinities of asymptotically hyperbolic manifolds.
\vspace{6pt}

\noindent
To introduce the fractional conformal Laplacian, we first recall some definitions in the theory of asymptotically hyperbolic metrics. Given \;$X=X^{n+1}$\; a smooth manifold with boundary $M=M^n$\; and \;$n\geq 2$ we say that \;$\varrho$\; is a defining function of the boundary \;$M$\; in \;
$X$,\; if 
$$
\varrho>0\;\;\text{ in }\;\;X,\;\;\varrho=0\;\;\text{ on }\;\;M\;\;\text{ and }\;\;d\varrho\neq 0\;\;\text{ on }\;\; M.
$$
A Riemannian metric \;$g^+$\; on \;$X$\; is said to be conformally compact, if for some defining function \;$\varrho$, the Riemannian metric 
\begin{equation}\label{eq:cmetric}
 g:=\varrho^2g^+
\end{equation}
 extends to \;$\ov{X}:=X\cup M$\; so that \;$(\ov{X}, \, g)$\; is a compact Riemannian manifold with boundary \;$M$ and interior \;$X$. Clearly this induces a conformal class of Riemannian metrics 
$$ [h]=[ g|_{TM}]$$ on \;$M$, where $TM$ denotes the tangent bundle of \;$M$, when the defining functions \;$\varrho$\; vary and the resulting conformal manifold \;$(M, [h])$\; is called conformal infinity of \;$(X, \;g^+)$. Moreover  a Riemannian metric \;$g^+$\; in \;$X$\; is said to be asymptotically hyperbolic, if it is conformally compact and its sectional curvature tends to \;$-1$\; as one approaches the conformal infinity of \;$(X, \;g^+)$,
 which is equivalent to 
\begin{equation*}
|d\varrho|_{\bar g}=1
\end{equation*}
 on\;$M$, see \cite{mazzeo1}, and in such a case \;$(X, \;g^+)$\; is called an asymptotically hyperbolic manifold. Furthermore a Riemannian metric \;$g^+$\; on \;$X$\; is said to be conformally compact Einstein or Poincar\'e-Einstein (PE), if it is asymptotically hyperbolic and satisfies the Einstein equation 
\begin{equation*}
Ric_{g^+}=-ng^+,
\end{equation*}
where $Ric_{g^+}$\; denotes the Ricci tensor of \;$(X, \;g^+)$.
\vspace{6pt}

\noindent
On one hand for every asymptotically hyperbolic manifold \;$(X, \;g^+)$\; and every choice of the representative \;$h$\; of its conformal infinity \;$(M, [h])$, there exists a  geodesic defining function \;$ y $\; of \;$M$\; in \;$X$\; such that  in a tubular neighborhood of \;$M$\; in \;$X$,  the Riemannian metric \;$g^+$\; takes the following normal form
\begin{equation}\label{eq:uniqdef}
g^+=\frac{d y ^2+h_{ y }}{ y ^2},
\end{equation}
where \;$h_{ y }$\; is a family of Riemannian metrics on \;$M$\; satisfying \;$h_0=h$\; and \;$ y $\; is the unique such a one in a tubular neighborhood of $M$.  Furthermore  we say that the conformal infinity \;$(M, \;[\hat h])$\; of an asymptotically hyperbolic manifold \;$(X, \;g^+)$\; is locally flat, if \;$h$\; is locally conformally flat, and clearly this is independent of the representative \;$h$\; of \;$[h]$. Moreover  we say that \;$(M, [h])$\; is umbilic, if \;$(M, h)$\; is umbilic in \;$(X, \; g)$ \;where \;$g$\; is given by \eqref{eq:cmetric} and \;$ y $\; is the unique geodesic defining function given by \eqref{eq:uniqdef}, and this is clearly independent of the representative \;$h$\; of \;$[h]$, as easily seen from the uniqueness of the normal form \eqref{eq:uniqdef} or Lemma 2.3 in \cite{gq}. Similarly we say that \;$(M, [h])$\; is minimal if \;$H_{g}=0$\; with \;$H_{g}$\; denoting the mean curvature of \;$(M,  \;h)$\; in \;$(\ov X, \;g)$ with respect to the inward direction, and this is again clearly independent of the representative of \;$h$\; of \;$[h]$, as easily seen from Lemma 2.3 in \cite{gq}. Finally we say that \;$(M, [h])$\; is totally geodesic, if \;$(M, [h])$\; is umbilic and minimal.                                                                                                                                         
\begin{rem}\label{eq:minimal} 
We remark that in the conformally compact Einstein case, \;$h_{ y }$\; as in \eqref{eq:uniqdef} has an asymptotic expansion which contains only even powers of \;$ y $, at least up to order \;$n$, see \cite{cg}.
In particular the conformal infinity \;$(M, [h])$\; of  any Poincar\'e-Einstein manifold \;$(X, g^+)$\; is totally geodesic.                                                                                                                                                                                                                                                                                                                                                                                                                                                                                                                                                                                                                                                                                                                                                                                                                                                                                                                                                                                                                                                                                                                                                                                                                                                                                                                                                                                                                                                                                                                                                                                                                                                                                                                                                                                                                                                                                                                                                                                                                                                                                                                                                                                                                                                                                                                                                                                                                                                                                                                                                                                                                                                                                                                                                                                                                                                                                                                                                                                                                                                                          \end{rem}
\begin{rem}
As every \;$2$-dimensional Riemannian manifold is locally conformally flat, we will say locally flat conformal infinity of a Poincar\'e-Einstein manifold  to mean just the conformal infinity of a Poincar\'e-Einstein manifold when the dimension is either \;$2$\; or which is further locally flat if the dimension is bigger than \;$2$. 
\end{rem}

\noindent
On the other hand, for any asymptotically hyperbolic manifold \;$(X, g^+)$\; with conformal  infinity \;$(M, [h])$, Graham-Zworsky\cite{gz}\; have attached a family of scattering operators \;$S(s)$\; which is a meromorphic family of pseudo-differential operators on \;$M$ defined on \;$\mathbb{C}$, by considering Dirichlet-to-Neumann operators for the scattering problem for \;$(X, \;g^+)$\; and a meromorphic continuation argument. Indeed  it follows from  \cite{gz} and \cite{mazmel} that for every \;$f\in C^{\infty}(M)$, and for every \;$s\in \mathbb{C}$\; such that \;$Re(s)>\frac{n}{2}$\;and\; \;$s(n-s)$ is not an \;$L^2$-eigenvalue of \;$-\D_{g^+}$, the following generalized eigenvalue problem
\begin{equation}\label{eq:geneig}
 -\D_{g^+}u-s(n-s)u=0\,\;\;\text{ in }\;\;X
\end{equation}
has a solution of the form
$$
u=F y ^{n-s}+G y ^s, \;\;\; F,\;G\in C^{\infty}(\ov X), \;\;\;F|_{ y =0}=f,
$$
where \;$ y $\; is given by  \eqref{eq:uniqdef} and for those values of \;$s$\; the scattering operator \;$S(s)$\; on \;$M$\; is defined  as
\begin{equation}\label{eq:dtnscat}
S(s)f=G|_{M}.
\end{equation}
Furthermore  using a meromorphic continuation argument, Graham-Zworsky\cite{gz} extend \;$S(s)$\; defined by \eqref{eq:dtnscat} to a meromorphic family of pseudo-differential operators on \;$M$\; defined on all \;$\mathbb{C}$\; and still denoted by $S(s)$\;with only a discrete set of poles including the trivial ones \;$s=\frac{n}{2}, \frac{n}{2}+1, \cdots, $\; which are simple poles of finite rank, and possibly some others corresponding to the \;$L^2$-eigenvalues of \;$-\D_{g^+}$. Using the regular part of the scattering operators \;$S(s)$, to any \;$\gamma\in (0, 1)$\; such that 
\begin{equation*}
\left(\frac{n}{2}\right)^2-\gamma^2<\l_1(-\D_{g^+})
\end{equation*}
with \;$\l_1(-\D_{g^+})$\: denoting the first eigenvalue of \;$-\D_{g_+}$,\; Chang-Gonzalez\cite{cg} have attached the following fractional order pseudo-differential operators referred to as fractional conformal Laplacians or fractional Paneitz operators
\begin{equation}\label{eq:fracopd}
P^{\gamma}[g^+, \;h]:=-d_\gamma S\left(\frac{n}{2}+\gamma\right),
\end{equation}
where \;$d_\gamma$\; is a positive constant depending only on \;$\gamma$\; and chosen  such that the principal symbol of \;$P^{\gamma}[g^+,  h]$\; is exactly the same as the one of the fractional Laplacian \;$(-\D_{h})^{\gamma}$, when
\begin{equation*}
X=\R^{n+1}_{+}, \;M=\R^{n},\;h=g_{\R^{n}} \;\;\text{ and }\;\; g^{+}=g_{\mathbb{H}^{n+1}}.
\end{equation*}
When there is no possible confusion with the  metric \;$g^+$, we just use the simple notation
\begin{equation*}
P^{\gamma}_{h}:=P^{\gamma}[g^+,  \;h].
\end{equation*}
Similarly to the other well studied conformally covariant differential operators, Chang-Gonzalez\cite{cg} associate to each \;$P^{\gamma}_{h}$\; the curvature quantity
\begin{equation*}
 Q^{\gamma}_{h}:=P^{\gamma}_{h}(1).
\end{equation*}
The \;$Q^\gamma_{h}$\; are referred to as fractional scalar curvatures, fractional \;$Q$-curvatures or simply \;$Q^\gamma$-curvatures. Of particular importance to conformal geometry is the  conformal covariance property verified by \;$P^{\gamma}_{h}$
\begin{equation}\label{eq:confinv}
P^{\gamma}_{h_u}(v)=v^{-\frac{n+2\gamma}{n-2\gamma}}P^{\gamma}_{h}(uv)
\;\;\text{ for }\;\; h_v=v^{\frac{4}{n-2\gamma}}  \;\;\text{ and }\;\; 0<v\in C^{\infty}(M).
\end{equation}
The fractional Yamabe problem is the problem of finding conformal metrics of  with constant \;$Q^\gamma$-curvature. As in the classical Yamabe problem, see \cite{sc}, its study deeply depends on the existence, regularity and sharp asymptotic of the Green's function of \;$P^{\gamma}_{h}$. 
\vspace{10pt}

\noindent
In this paper, we show existence, regularity and sharp asymptotics of the Poisson kernel \;$K_g$\; and Green's functions  \;$\Gamma_g$\; under weighted Neumann boundary conditions of the Chang-Gonzalez\cite{cg} extension problem associated to \;$P^{\gamma}_{h}$ and the Green's function\;$G_h$\; of \;$P^{\gamma}_{h}$. Indeed we prove:
\begin{thm}
\label{Greens_function_asymptoticsgamma1}$_{}$\\
Let \;$(X,  \;g^{+})$\; be an asymptotically hyperbolic manifold with 
conformal infinity  \;$(M, [h])$\; of dimension \;\;$n\geq 2$. If
\begin{equation*}
\frac{1}{2}\neq \gamma \in (0,1) \;\;\text{ and } \;\,\;
\lambda_{1}(-\Delta_{g^{+}})>s(n-s)\;\;
\text{ for } s=\frac{n}{2}+\gamma,
\end{equation*}
then the Poison kernel \;$K_{g}$\; and the Green's functions \;$\Gamma_{g}$\; and \;$G_{h}$\;  
respectively for
\begin{equation*}
\begin{cases}
D_{g}U=0 \;\;\text{ in }\;\; X \\
U=f \,\,\;\text{ on } \;\;M
\end{cases} 
\quad\quad
\begin{cases}
D_{g}U=0\;\; \text{ in } \;\;X \\
-d_{\gamma}^*\lim_{y\rightarrow 0}y^{1-2\gamma}\partial_{y}U=f \;\;\text{ on }\;\; M
\end{cases} 
\quad \text{and}\quad
\begin{cases}
P_h^\gamma u=f\;\;\text{ on }\;\;M
\end{cases}
\end{equation*}
exist and we may expand in \;$g$-normal 
Fermi-coordinates 
around $\xi \in M$
\begin{enumerate}[label=(\roman*)]

 \item 
 \quad
 $
K_{g}(z,\xi)
\;\in
\eta_{\xi}(z)
\left( p_{n, \gamma}
\frac{y^{2\gamma}}{\vert z 
\vert^{n+2\gamma}}+\sum^{2m+5-2\gamma}_{l=-n-2\gamma}y^{2\gamma}H_{1+l}(z)
\right)
+
y^{2\gamma}C^{2m,\alpha}(X)
$
\item
 \quad 
$
\Gamma_{g}(z,\xi)
\;\in
\eta_{\xi}(z)
\left( 
\frac{g_{n, \gamma}}{\vert z \vert^{n-2\gamma}}+\sum^{2m+3}_{l=-n}H_{1+2\gamma+l}(z)
\right)
+
C^{2m,\alpha}(X)
$

 \item \quad 
$
G_{h}(x,\xi)
\;\in
\eta_{\xi}(x)
\left( 
\frac{g_{n, \gamma}}{\vert x \vert^{n-2\gamma}}+\sum^{2m+3}_{l=-n}H_{1+2\gamma+l}(x)
\right)
+
C^{2m,\alpha}(M)
$

\end{enumerate}
with \;$H_{l}\in C^{\infty}(\R^{n+1}_+\setminus \{0\})$\; being homogeneous of order 
\;$l$,   \;$\eta_{\xi}$\; as in \eqref{etaxi}, \;$p_{n, \gamma}$\; is as in \eqref{pngamma}, and \;$g_{n, \gamma}$\; is as in \eqref{gngamma}, provided \;$H_{g}=0$.
\end{thm}
\vspace{6pt}

\noindent
In the case of  locally flat conformal infinities of Poincare-Einstein manifolds, we have:
\begin{thm}
\label{cor_kernels_for_poincare_einstein_metrics1}$_{}$\\
\noindent
Let \;$(X,  \;g^{+})$\; be a Poincar\'e-Einstein manifold with 
conformal infinity  \;$(M, [h])$\; of dimension \;$n= 2$\; or \;$n\geq 3$\; and \;$(M, [h])$\; is locally flat. If
\begin{equation*}
\frac{1}{2}\neq \gamma \in (0,1) \;\;\text{ and } \;\,\;
\lambda_{1}(-\Delta_{g^{+}})>s(n-s)\;\;
\text{ for } s=\frac{n}{2}+\gamma,
\end{equation*}
then the Poisson kernel \;$K_{g}$\; and the Green's functions \;$\Gamma_{g}$\; and \;$G_{h}$\;  
respectively for
\begin{equation*}
\begin{cases}
D_{g}U=0 \;\;\text{ in }\;\; X \\
U=f \,\,\;\text{ on } \;\;M
\end{cases} 
\quad\quad
\begin{cases}
D_{g}U=0\;\; \text{ in } \;\;X \\
-d_{\gamma}^*\lim_{y\rightarrow 0}y^{1-2\gamma}\partial_{y}U=f \;\;\text{ on }\;\; M
\end{cases} 
\quad \text{and}\quad
\begin{cases}
P_h^\gamma u=f\;\;\text{ on }\;\;M
\end{cases}
\end{equation*}
are respectively of class \;$y^{2\gamma}C^{2,\alpha}$\; and \;$C^{2,\alpha}$\; away from the 
singularity and admit for every \;$a\in M$ locally in \;$g_{a}$-normal 
Fermi-coordinates an expansion around $a$
\begin{enumerate}[label=(\roman*)]
 \item \quad 
$
K_{a}(z)\;
\in
p_{n, \gamma}\frac{y^{2\gamma}}{\vert z \vert^{n+2\gamma}}
+
y^{2\gamma}H_{-2\gamma}(z)
+
y^{2\gamma}H_{1-2\gamma}(z)
+
y^{2\gamma}H_{2-2\gamma}(z)
+
y^{2\gamma}C^{2,\alpha}(X)
$

 \item \quad 
$
\Gamma_{a}(z)\;
\in
\frac{g_{n, \gamma}}{\vert z \vert^{n-2\gamma}}+H_{2\gamma}(z)+H_{1+2\gamma}(z)
+
C^{2,\alpha}(X)
$

 \item \quad 
$
G_{a}(x)\;
\in
\frac{g_{n, \gamma}}{\vert x \vert^{n-2\gamma}}+H_{2\gamma}(x)+H_{1+2\gamma}(x)
+
C^{2,\alpha}(M),
$

\end{enumerate}
where \;$g_a$\;is as in  \eqref{dfmga},  \;$K_{a}=K_{g_{a}}(\cdot, a)$, \;$\Gamma_{a}=\Gamma_{g_{a}}(\cdot, 
a)$ and \;$G_a=G_{h_a}(\cdot,a)$\; and \;$H_{k}\in C^{\infty}(\overline{\R^{n}_{+}}\setminus \{0\})$\; are 
homogeneous of degree \;$k$. 
\end{thm}
\vspace{10pt}

\noindent
To prove  Theorem \ref{Greens_function_asymptoticsgamma1} and Theorem \ref{cor_kernels_for_poincare_einstein_metrics1}, we use the method of Lee-Parker\cite{lp} of killing deficits successively. However difficulties arise due the the rigidity involved in the problem (see \eqref{eq:uniqdef}) and the lack of classical regularity theory. To overcome the rigidity issue, we work with the space of homogeneous functions rather than the one of polynomials as done in \cite{lp}. To handle the regularity issue, we show some higher order regularity results for the Dirichlet problem and the weighted Neumann boundary problem of the Chang-Gonzalez\cite{cg} extension problem for \;$P_h^{\gamma}$\; which are of independent interest, see Proposition \ref{poisson_regularity} and Proposition \ref{green_regularity}. We point out that even if the estimates in Proposition \ref{poisson_regularity} and Proposition \ref{green_regularity} are weak, they are enough for our purpose and in turn get improved by the estimates of the Poisson kernel and Green's function in Theorem \ref{Greens_function_asymptoticsgamma1} and Theorem \ref{cor_kernels_for_poincare_einstein_metrics1} that they imply. On the other hand, we would like to emphasize that (ii) of Theorem \ref{cor_kernels_for_poincare_einstein_metrics1} answers the first part of the Conjecture of Kim-Musso-Wei\cite{kmw1} about the asymptotics of \;$\Gamma_a$\; and gives the definition of the fractional mass, see our work \cite{martndia3}, Definition 4.3 and Lemma 4.1.

\vspace{10pt}

\noindent
The structure of the paper is as follows:   In Section \ref{notation_and_prelimiaries} we fix some notations. In Section \ref{nonhom} we develop a non-homogeneous extension of some aspects of the works of Chang-Gonzalez\cite{cg} and Graham-Zworsky\cite{gz}. It is divided in two subsections. In the first one, namely Subsection \ref{eq:nhscatdeg}, we develop a non-homogeneous scattering theory, define the associated non-homogeneous fractional operator and its relation to a non-homogeneous uniformly degenerate boundary value problem. In Subsection \ref{confnonhom} we discuss the conformal property of the non-homogeneous fractional operator. We point out that Section \ref{nonhom} even being of independent interest contains estimates which are used in Section \ref{pceinsten} and in \cite{martndia3}, and a regularity result that we use in \cite{martndia3} . Section \ref{fundso} is concerned with the study of the Poisson kernel \;$K_g$  and  the Green's function \;$\Gamma_g$\; under weighted Neumann boundary conditions of the Chang-Gonzalez extension problem \;of $P_h^\gamma$, and the Green's function \;$G_h$\; of $P_h^\gamma$\; all in the general case of asymptotically hyperbolic manifolds with minimal conformal infinity. In Section \ref{pceinsten} we sharpen the results obtained in Section \ref{fundso} in the particular case of a locally flat conformal infinity of a Poincar\'e-Einstein manifold.

\section{Notations and preliminaries }\label{notation_and_prelimiaries} 
In this section we fix some notations. First of all let \;$X=X^{n+1}$\; be a manifold of dimension \;$n+1$\; with boundary \;$M=M^{n}$\; and closure \;$\ov{X}$\; with\;$n\geq 2$. 
\vspace{10pt}

\noindent
In the following, for any Riemannian metric \;$\bar h$\; defined on \;$M$, $a\in M$\; and \;$r>0$, we use the notation \;$B^{\bar h}_{r}(a)$\; to denote the geodesic ball with respect to $\bar h$\; of radius \;$r$\;and center \;$a$.  We also denote by \;$d_{\bar h}(x,y)$\; the geodesic distance with respect to \;$\bar h$\; between two points \;$x$\;and \;$y$\; of \;$M$.  $inj_{\bar h}(M)$\;stands for the injectivity radius of \;$(M, \bar h)$. $dV_{\bar h}$\;denotes the Riemannian measure associated to the metric\;$\bar h$\; on \;$M$. For \;$a\in M$ we use the notation \;$\exp^a_{\bar h}$\; to denote the exponential map with respect to \;$\bar h$\; on \;$M$. 
\vspace{6pt}

\noindent
Similarly for any Riemannian metric \;$\bar g$\; defined on \;$\ov{X}$, $a\in M$\; and \;$r>0$ we use the notation \;$B^{\bar g, + }_{r}(a)$\; to denote the geodesic half ball with respect to \;$\bar g$\; of radius \;$r$\;and center \;$a$.  We also denote by \;$d_{\bar g}(x,y)$\; the geodesic distance with respect to \;$\bar g$\; between two points \;$x\in M$\;and \;$y\in \ov {X}$.  $inj_{\bar g }(\ov{X})$\;stands for the injectivity radius of \;$(\ov{X}, \bar g)$. $dV_{\bar g}$\;denotes the Riemannian measure associated to the metric\;$\bar g$\; on \;$\ov{X}$. For \;$a\in M^{n}$\; we use the notation \;$\exp_a^{\bar g, +}$\; to denote the exponential map with respect to \;$\bar g$\; on \;$\ov{X}$. 
\vspace{6pt}

\noindent
$\N$\;denotes the set of nonnegative integers, $\N^*$\; the set of positive integers and  for $k\in \N^*$, $\R^k$\;stands for the standard $k$-dimensional Euclidean space,  $\R^k_+$ the open positive half-space of $\R^k$, and $\bar \R^k_+$ its closure in $\R^k$. For simplicity we use the notation \;$\R_+:=\R^1_+$, and $\bar \R_+:=\bar \R^1_+$. For $r>0$
we denote respectively
\begin{equation*}
B^{\R^k}_r(0) \;\; \text{ and }\;\; B^{\R^k_+}_r(0)=B_r^{\R^{k}}(0)\cap  \R^k_+\simeq ]0, r[\times B_r^{\R^{k-1}}(0)
\end{equation*}
the open and open upper half ball of \;$\R^k$\; of center \;$0$\; and radius \;$r$, and set
 \;$B_r=B_r^{\R^n}$\; and $B_r^+=B_r^{\R^{n+1}_+}$. 
\vspace{6pt}

\noindent
For \;$p\in \N^*$, let \;$M^p$\; denotes the Cartesian product of \;$p$\; copies of \;$M$. We define \;$(M^2)^*:=M^2\setminus Diag(M^2)$, where \;$Diag(M^2)=\{(a, a): \,a\in M\}$\; is the diagonal of \;$M$.
\vspace{6pt}

\noindent
For \;$1\leq p\leq \infty,\;k\in \N$, $s\in \R_+$, $\beta\in  ]0, 1[$ and \;$\bar h$\; a Riemannian metric defined on \;$M$, \begin{equation*}
L^p(M, \bar h), \;W^{s, p}(M, \bar h), \;C^k(M, \bar h)\;\;\text{ and }\;\;C^{k, \beta}(M, \bar h)
\end{equation*}
stand respectively for the standard $p$-Lebesgue and \;$(s, p)$-Sobolev space, $k$-continuously differentiable space and \;$k$-continuously differential space of H\"older exponent \;$\beta$, all on \;$M$ and with respect to \;$\bar h$, if the definition required a metric structure.  Similarly for\;$1\leq p\leq \infty,\;k\in \N$, $s\in \R_+$, $\beta\in  ]0, 1[$ and \;$\bar g$\; a Riemannian metric defined on \;$\ov{X}$, 
\begin{equation*}
L^p_{f}(\ov{X}, \bar g), \;W^{s, p}_{f}(\ov{X}, \bar g ), \;C^k(\ov{X}, \bar g)\text{ and }\;C^{k, \beta}(\ov{X}, \bar g) 
\end{equation*}
stand respectively for the weighted \;$p$-Lebesgue and \;$(s, p)$-Sobolev space, \;$k$-continuously differentiable space and \;$k$-continuously differential space of H\"older exponent \;$\beta$, all on \;$\ov{X}$, and as above with respect to \;$\bar g$\; and a measurable function \;$f>0$\; on \;$X$\;, if required. For precise definitions and properties see \cite{aubin}, \cite{dinpalval}, \cite{gold}, \cite{gt} and \cite{s}.
$C_0^{\infty}(X)$ means element in \;$C^{\infty}(X)$\; vanishing on \;$M$\; to infinite order.
\vspace{6pt} 

\noindent
For \;$\epsilon>0$\;  and small \;$o_{\epsilon}(1)$\; means quantities which tend to \;$0$\; as \;$\epsilon$\; tends to \;$0$. $O(1)$ stands for quantities which are bounded.  For \;$x\in \R$\; we use the notation \;$O(x)$\; and 
\;$o_{\epsilon}(x)$  to mean  respectively \;$|x|O(1)$\; and \; $|x|o_{\epsilon}(1)$.                                                                                                                                                                                                                                                                                                                                                                                                                                                                                                                                                                                                                                                                                                                                                                                                                                                                                                                                                                                                                                                                                                    Large positive constants are usually denoted by \;$C$\; and the value of \;$C$\; is allowed to vary from formula to formula and also within the same line. Similarly small positive constants are denoted by \;$c$\; and their values may vary from formula to formula and also within the same line. 
\vspace{6pt}

\noindent
We define
\begin{equation}\label{dsgamma}
d_{\gamma}^*=\frac{d_\gamma}{2\gamma}, 
\end{equation}
cf. \eqref{eq:fracopd}. Furthermore, we set
\begin{equation}\label{defc3}
c_{n, 3 }^{\gamma}=\int_{\R^{n}}\left(\frac{1}{1+|x|^2}\right)^{\frac{n+2\gamma}{2}}dx,\;\;
\end{equation}
and
\begin{equation}\label{pngamma}
p_{n,\gamma}=\frac{1}{c_{n,3}^{\gamma}}
\end{equation}
\vspace{4pt}

\noindent
Let\;$(X, g^+)$  be an asymptotically hyperbolic manifold of dimension \;$n+1$\; with  \;$n\geq 2$\; and minimal conformal infinity \;$(M, [h])$. Then, because of \eqref{eq:uniqdef} and minimality of the conformal infinity, we can consider a geodesic defining function \;$y$\; splitting the metric 
\begin{equation*}
g=y^{2}g^{+}, \;\; g=dy^{2}+h_{y}\;\;\text{near}\;\;M\;\;\text{ and } \;\;h=h_{y}\lfloor_{M}
\end{equation*} 
in such a way, that \;$H_{g}=0$.
Moreover  using the existence of conformal normal coordinates, cf.  \cite{gun}, there exists for every \;$a\in M$\; a conformal factor
\begin{equation}\label{conformal_factor_properties}
0<u_{a}\in C^{\infty}(M)\;\;\text{ satisfying }\;\;\frac{1}{C}\leq u_{a}\leq 
C, \;\;u_{a}(a)=1\;\;\text{ and }\;\;\nabla u_{a}(a)=0,
\end{equation}
inducing a conformal normal coordinate system close to \;$a$\; on \;$M$, in 
particular in normal coordinates with respect to  
$$h_{a}=u_{a}^{\frac{4}{n-2\gamma}}h$$  we have for some small \;$\epsilon>0$
\begin{equation*}
h_{a}=\delta+O(\vert x \vert^{2}), \;\;\det h_{a}\equiv 1\;\;\text{ on } \;\;
B_{\epsilon}^{h_a}(a).
\end{equation*}
As clarified in 
Subsection \ref{confnonhom} the conformal factor \;$u_{a}\;$ then 
naturally extends onto \;$X$ via 
$$u_{a}=(\frac{y_{a}}{y})^{\frac{n-2\gamma}{2}},$$
where \;$y_{a}$\; close to the 
boundary \;$M$ is the unique geodesic defining function, for which
\begin{equation*}
g_{a}=y_{a}^{2}g^{+}, \;\; g_{a}=dy_{a}^{2}+h_{a, y_a}\;\;\text{near}\;\;M\;\;\text{ with } 
\;\;h_{a}=h_{a, y_a}\lfloor_{M}
\end{equation*} 
and there still holds \;$H_{g_{a}}=0$. Consequently
\begin{equation*}
g_{a}=\delta+O(y+\vert x \vert^{2})\;\;\text{ and }\;\; \det g_{a}=1+O(y^{2}) \;\;\text{in 
} \;\;B_{\epsilon}^{g_a, +}(a).
\end{equation*}

\section{Non-homogeneous scattering theory}\label{nonhom}
In this section we extend some aspects of the works of Chang-Gonzalez\cite{cg} and Graham-Zworsky\cite{gz}  to a non-homogeneous setting and in the general framework of asymptotically  hyperbolic manifolds. It is of independent interest, but in it we derive estimates that are used in Section \ref{pceinsten} and \cite{martndia3}, and an existence and regularity result used in \cite{martndia3} to construct barrier solutions in order to compare different types of bubbles via maximum principle.  We divide this section in two subsections.
\subsection{Scattering operators and uniformly degenerate equations}\label{eq:nhscatdeg}
In this subsection we extend some parts of the works of Chang-Gonzalez\cite{cg} and Graham-Zworski\cite{gz} to a non-homogeneous setting in the context of asymptotically hyperbolic manifolds. First of all  let \;$(X, g^{+})$\; be an asymptotically hyperbolic manifold with conformal infinity \;$(M, [h])$\; and \;$y$\; the unique geodesic defining function associated to \;$h$\; given by \eqref{eq:uniqdef}. Then we have the normal form 
\begin{equation*}\begin{split}
\;y^{2}g^{+}=g=dy^{2}+h_{y} \;\;\text{ near }\;\;M
\end{split}\end{equation*}
with
\;$
y>0\;\; \text{ in }  \;\;X,\;\;y=0\;\;\text{ on }\;\; M\;\; \text{ and }\;\;\vert dy \vert_{g}=1 \;\;\text{ near }\;\; M.
$ Furthermore  let
\begin{equation*}
\square_{g^{+}}=-\Delta_{g^{+}}-s(n-s),
\end{equation*}
where by definition
\begin{equation*}
s=\frac{n}{2}+\gamma,\;\;\gamma\in (0,1), \;\,\gamma\neq \frac{1}{2}\;\; \text{and}\; \;\;s(n-s)\in (0,\frac{n^{2}}{4}).
\end{equation*} 
According to Mazzeo and Melrose \cite{mazzeo1}, \cite{mazzeo2}, \cite{mazmel}
\begin{equation*}
\sigma(-\Delta_{g^{+}})=\sigma_{pp}(-\Delta_{g^{+}})\cup [\frac{n^{2}}{4}, \infty),
\; \sigma_{pp}(-\Delta_{g^{+}})\subset (0,\frac{n^{2}}{4}),
\end{equation*}
where \;$\sigma(-\Delta_{g^{+}})$\; and \;$\sigma_{pp}(-\Delta_{g^{+}})$ are respectively the spectrum and  the pure point spectrum of \;$L^{2}$-eigenvalues of \;$-\Delta_{g^{+}}$.  
Using the work of  Graham-Zworski\cite{gz}, see equation (3.9) therein, we may solve
\begin{equation*}\begin{split}
\begin{cases}
\square_{g^{+}}u=f\;\;\text{ in }\;\;X \\
y^{s-n}u=\underline{v} \;\;\text{ on } \;\;M
\end{cases}
\end{split}\end{equation*}
for \;$s(n-s)\not \in \sigma_{pp}(-\Delta_{g^{+}})$\; and \;$f\in y^{n-s+1}C^{\infty}(\ov X)+y^{s+1}C^{\infty}(\ov X)$\; in the form
\begin{equation*}\begin{split}
\begin{cases}
u=y^{n-s}A+y^{s}B \;\;\text{ in }\;\;X\\
A,\; B\in C^{\infty}(\ov X),\;\;\,A=\underline{v} \;\;\text{ on } \;\;M.
\end{cases}
\end{split}\end{equation*}
As in the case \;$f=0$, which corresponds to the generalized eigenvalue problem of Graham-Zworsky\cite{gz}, this gives rise to a Dirichlet-to-Neumann map \;$S_{f}(s)$\; via
\begin{equation*}\begin{split}
\underline{v}=A\lfloor _{M}\longrightarrow\-B\lfloor _{M}=\overline{v},
\end{split}\end{equation*}
which we refer to as non-homogeneous scattering operator and denote it by \;$S_{f}(s)$. Clearly $S_0(s)=S(s)$ and $S_{f}(s)$ is invertible, since the standard scattering operator \;$S_0(s)$\; 
is invertible, cf. equation (1.2) in \cite{jb}. We define the non-homogeneous fractional operators by
$$
P^{\gamma}_{f,h}=-d_\gamma S_{f}(s),
$$
where \;$d_\gamma$\; is as in \eqref{eq:fracopd}. Following \cite{gq} we find by conformal covariance of the conformal Laplacian that 
\begin{equation}\begin{split}\label{transformation_scattering_extension}
\square_{g^{+}}u=f
\xLeftrightarrow{U=y^{s-n}u} D_{g}U=y^{-s-1}f,
\end{split}\end{equation}
where
\begin{equation}\label{eqdg}
\begin{split}
D_{g}U=-div_{g}(y^{1-2\gamma}\nabla_{g}U)+E_{g}U
\end{split}
\end{equation}
and with \;$L_g=-\Delta_{g}+\frac{R_{g}}{c_{n}}$\; denoting the conformal Laplacian on \;$(X, g)$
\begin{equation}\begin{split}\label{Eg_general}
E_{g} :=y^{\frac{1-2\gamma}{2}}L_{g}y^{\frac{1-2\gamma}{2}} -(\frac{R_{g^{+}}}{c_{n}}+s(n-s))y^{(1-2\gamma)-2}
,\;\;
c_{n}=\frac{4n}{n-1}.
\end{split}
\end{equation}
Thus we find for \;$\phi,\psi \in C^{\infty}(\ov X)$, that 
\begin{equation*}\begin{split}
\begin{cases}
\square_{g^{+}}u=y^{n-s+1}\phi +y^{s+1}\psi\;\;\text{ in }\;\;X \\
y^{s-n}u=\underline{v} \;\; \text{ on } \;\;M
\end{cases}
\xLeftrightarrow{U=y^{s-n}u}
\begin{cases}
D_{g}U=y^{-2\gamma}\phi + \psi\;\;\text{ in }\;\;X \\
U=\underline{v} \;\;\text{ on } \;\;M
\end{cases}.
\end{split}\end{equation*}
Note, that such a solution \;$U$\; is of the form
\begin{equation*}
\begin{split}
U=A+By^{2\gamma}=\sum A_{i}y^{i}+\sum B_{i}y^{i+2\gamma}+U_0 
\end{split}
\end{equation*} 
for some \;$U_0\in C^{\infty}_0(X)$\; and has  principal terms
\begin{equation*}
\begin{cases} 
\underline{v}+\overline{v}  y^{2\gamma}\;\; \text{ for }\;\; \gamma<\frac{1}{2}\\
\underline{v}+A_{1}y+\overline{v} y^{2\gamma}\;\; \text{ for }\;\; \gamma>\frac{1}{2}.
\end{cases}
\end{equation*}
As for the case $\gamma>\frac{1}{2}$, expanding the boundary metric  
\;$h_y$,\;  we find 
\begin{equation*}
\begin{split}
h_{y}=h_{0}+h_{1}y+O(y^{2}) \;\;\text{ with }\;\; h_{1}=2\Pi_g
\end{split}
\end{equation*} 
and \;$\Pi_g$\; denoting the second fundamental form of \;$(M, h)$\; in $(\ov{X}, g)$. Still according to \cite{gz} we may solve
\begin{equation*}\begin{split}
\begin{cases}
\square_{g^{+}}u=\square_{g^{+}}u=y^{n-s+2}\phi +y^{s+1}\psi\;\;\text{ in }\;\;X \\
y^{s-n}u=\underline{v} \;\;\text{on} \;\;M
\end{cases}
\end{split}\end{equation*}
for \;$\phi,\psi \in C^{\infty}(\ov X)$ in the form
\begin{equation*}\begin{split}
\begin{cases}
u=y^{n-s}A+y^{s}B \;\;\text{ in }\;\;X\\
A, \;B\in C^{\infty}(\ov X),\;\;A=\underline{v} \;\;\text{ on } \;\;M
\end{cases}
\end{split}\end{equation*}
with asymptotic 
\begin{equation*}
\begin{split}
A=\sum A_{i}y^{i}, \;\; A_{0}=\underline{v}, \;\;A_{1}=0
\end{split}
\end{equation*} 
at a point, where \;$H_g=0$, i.e. the mean curvature vanishes. Thus for \;$\gamma>\frac{1}{2}$
\begin{equation*}\begin{split}
\begin{cases}
\square_{g^{+}}u=y^{n-s+2}\phi +y^{s+1}\psi\;\;\text{ in }\;\;X \\
y^{s-n}u=\underline{v}  \;\;\text{ on } \;\; M
\end{cases}
\xLeftrightarrow{U=y^{s-n}u}
\begin{cases}
D_{g}U=y^{1-2\gamma}\phi+\psi\;\;\text{ in }\;\;X \\
U= \underline{v} \;\;\text{ on } \;\;M
\end{cases}
\end{split}\end{equation*}
with principal terms
\begin{equation*}
U=\underline{v}+\overline{v} y^{2\gamma}+o(y^{2\gamma})
\end{equation*}
at a point with \;$H_{g}=0$ - just like in the case \;$\gamma<\frac{1}{2}$ - and there holds
$
\overline{v}=\frac{1}{2\gamma}\lim_{y\to 0}y^{1-2\gamma}\partial_{y}U.
$
\vspace{6pt}

\noindent
We summarize the latter discussion in the following proposition.
 \begin{pro}\label{prop_scattering}
 Let \;$(X, g^{+})$\; be a \;$(n+1)$-dimensional asymptotically hyperbolic manifold with conformal infinity \;$(M, [h])$\; of dimension \;$n\geq 2$ \;being minimal in case \;$\gamma\in(\frac{1}{2},1)$\; and \;$y$ the unique geodesic defining function associated to \;$h$\; given by \eqref{eq:uniqdef}. Assuming that $$s=\frac{n}{2}+\gamma,\;\;\gamma\in (0,1), \;\,\gamma\neq \frac{1}{2}, \; \;\;s(n-s)\not \in \sigma_{pp}(\Delta_{g^{+}})$$
and \;$f\in y^{n-s+2}C^{\infty}(\ov {X})+y^{s+1}C^{\infty}(\ov {X})$, then for every \;$\underline{v}\in C^{\infty}(M)$
$$
P^{\gamma}_{f,h}(\underline{v})=-d_\gamma^*\lim_{y\to 0}y^{1-2\gamma}\partial_{y}U^{f},
$$
where \;$U^{f}$\; is the unique solution to
\begin{equation*}
\begin{cases}
D_{g}U=y^{-s-1}f\;\;\text{in}\;\;X
\\
U=\underline{v}\;\;\text{on} \;\;M
\end{cases}
\end{equation*} 
and \;$d_\gamma^*$ is as in \eqref{dsgamma}.  Moreover  $U^{f}$\; satisfies
 $$U^{f}=A+y^{2\gamma}B,\; \;\;A, \;B\in C^{\infty}(\ov{X})$$ and \;$A$\; and \,$B$\; satisfy the asymptotics
\begin{equation*}
\begin{cases}
A=\sum A_{i}y^{i}, \;\;\,A_{i}\in C^{\infty}(M),\;\; A_{0}=\underline{v}\; \text{ and }\;\; A_{1}=0 \\
B=\sum B_{i}y^{i}, \;\;B_{i}\in C^{\infty}(M)\; \;\;\text{ and }\;\; -d_\gamma B_{0}=-d_\gamma\overline{v}=P^{\gamma}_{f,h}(\underline{v}),
\end{cases}
\end{equation*}
where \;$d_\gamma$\; is as in \eqref{eq:fracopd}, hence
$U^{f}=\underline{v}+\overline{v} y^{2\gamma}+o(y^{2\gamma}).$
\end{pro}
\subsection{Conformal property of the non-homogeneous scattering operator}\label{confnonhom}
In this subsection we study the conformal property of the non-homogeneous scattering operator \;$P_{h, f}^{\gamma}$\; of the previous subsection.  To this end we first 
consider as background data \;$(X, g^{+})$\; with conformal  infinity \;$(M, [h])$\; with \;$n\geq 2$\; and \;$y$\; the associated unique geodesic definition function  such that 
\begin{equation*}
g=y^{2}g^{+}, \;\; g=dy^{2}+h_{y}\;\; \text{ close to \;$M$\; and }\;\;h=g\lfloor_{M}
\end{equation*}
as in \eqref{eq:uniqdef}. From \eqref{Eg_general} it is easy to see, that in \;$g$-normal Fermi coordinates \;$(y, x)$
\begin{equation}\label{Eg_local}
E_{g}=\frac{n-2\gamma}{2}\frac{\partial_{y}\sqrt{g}}{\sqrt{g}}y^{-2\gamma}\;\;
\text{ close to }\;\; M.
\end{equation}
We assume further that \;$(M, [h])$\; is minimal and \;$\square_{g^+}$\; is positive, i.e. 
\begin{equation*}
H_{g}=0 \;\; \text{ and }\;\; \lambda_{1}(-\Delta_{g^{+}})>s(n-s).
\end{equation*}
Then \;$\partial_{y}\sqrt{g}=0$\; on \;$M^{n+1}$\; and we may assume \;
\begin{equation}\label{sqrt_g_sim_yC_infty}
\partial_{y}\sqrt{g}\in 
yC^{\infty}(\ov {X}) 
\end{equation} whence 
$D_{g}$\; is well defined on 
$$
W^{1,2}_{y^{1-2\gamma}}=W^{1,2}_{y^{1-2\gamma}}(X,g)
=
\overline{C^{\infty}(X)}^{\Vert \cdot \Vert_{W^{1,2}_{y^{1-2\gamma}}(X,g)}},
\;
\Vert u \Vert_{W^{1,2}_{y^{1-2\gamma}}(X,g)}^{2}
=
\int_{X} y^{1-2\gamma}(\vert du \vert^{2}_{g}+u^{2})dV_{g}
$$
and becomes positive under Dirichlet condition, cf. \eqref{transformation_scattering_extension}, so
\begin{equation*}
\partial_{y}\sqrt{g}\in yC^{\infty}(\ov{X})\;\;\text{ and }\;\;\langle\cdot,\cdot \rangle_{D_{g}} \simeq \langle \cdot,\cdot \rangle_{W^{1,2}_{y^{1-2\gamma}}}.
\end{equation*} 
Let us consider now a conformal metric \;$\tilde h=\varphi^{\frac{4}{n-2\gamma}}h$\; on \;$M$. We then find a unique geodesic defining function \;$\tilde y>0$, precisely unique in a tubular neighborhood  of \;$M$, such that
\begin{equation*}
\tilde g=d\tilde y^{2}+\tilde h_{y}\;\;\text{ close to }\;\;M, \;\;\tilde{y}^{-2}\tilde g=g^{+}=y^{-2}g 
\;\;\text{ and }\;\;
\tilde h=\varphi^{\frac{4}{n-2}}h=(\frac{\tilde y}{y})^{2}h\;\;\text{ on }\;\;M.
\end{equation*} 
So we may naturally extend \;$\varphi=(\frac{\tilde y}{y})^{\frac{n-2\gamma}{2}}$\; onto \;$X$\; and
by the conformal  relation 
\begin{equation*}
\tilde g=(\frac{\tilde y}{y})^{2}g=\varphi^{\frac{4}{n-2\gamma}}g,
\end{equation*}
we still have \;$\langle \cdot,\cdot \rangle_{D_{\tilde g}}\simeq \langle \cdot,\cdot \rangle_{W^{1,2}_{\tilde y^{1-2\gamma}}}$. Putting \;$\tilde y=\alpha y$, the equation
\begin{equation*}
\vert d y\vert^{2}_{g}
=
1
=
\vert d \tilde y\vert^{2}_{\tilde g}
=
1
+
2\frac{y}{\alpha}\langle d\alpha,dy\rangle_{g}+ (\frac{y}{\alpha})^{2}\vert d\alpha \vert_{g}^{2}
\end{equation*} 
for the geodesic defining functions implies 
$\partial_{y}\alpha=-\frac{1}{2}\frac{y}{\alpha}\vert d\alpha\vert^{2}_{g}$. Since $\tilde g =\alpha^{2}g$ by definition, we firstly find
$H_{g}=0 \Longrightarrow H_{\tilde g}=0$, i.e. minimality is preserved  as already observed by Gonzalez-Qing\cite{gq}, and secondly  \;$\tilde y=\alpha_{0}y+O(y^{3})$.
Thus on the one hand side  the properties
\begin{equation*}
\partial_{\tilde y}\sqrt{\tilde g}\in \tilde y C^{\infty}
\text{ and }
\langle \cdot,\cdot \rangle_{D_{\tilde g}}\simeq \langle \cdot,\cdot \rangle_{W^{1,2}_{\tilde y^{1-2\gamma}}}\;\;
\end{equation*} 
are preserved under a conformal change of the metric on the boundary. Moreover we obtain a conformal transformation for the extension operators \;$D_{\tilde g}$\; and $\;D_{g}$\; subjected to Dirichlet and weighted Neumann boundary conditions. Put \;$\tilde u=(\frac{y}{\tilde y})^{n-s}u$. 
As for the  Dirichlet case, \eqref{transformation_scattering_extension} directly shows
\begin{equation*}
\begin{cases}
D_{g}u=f \;\;\text{ in } \;\;X \\
u=  v\;\;\text{ on } \;\;M
\end{cases}
\Longleftrightarrow
\begin{cases}
D_{\tilde g}\tilde u=(\frac{y}{\tilde y})^{s+1}f \;\;\text{ in } \;\;X \\
\tilde u=(\frac{y}{\tilde y})^{n-s}v \;\;\text{ on } \;\;M.
\end{cases}
\end{equation*} 
Moreover there holds
$$
\lim_{y\rightarrow 0}y^{1-2\gamma}\partial_{y}u=v
\Longleftrightarrow 
\lim_{\tilde y\rightarrow 0}\tilde y^{1-2\gamma}\partial_{\tilde y}\tilde u=(\frac{y}{\tilde y})^{n-s+2\gamma}v,
$$
since \;$\tilde y=\alpha_{0}y+O(y^{3})$, whence for the weighted Neumann case we obtain
\begin{equation*}
\begin{cases}
D_{g}u=f \;\;\text{ in } \;\;X \\
\lim_{y\rightarrow 0}y^{1-2\gamma}\partial_{y}u=v \;\;\text{ on } \;\;M
\end{cases}
\Longleftrightarrow
\begin{cases}
D_{\tilde g}\tilde u=(\frac{y}{\tilde y})^{s+1}f \;\;\text{ in } \;\;X \\
\lim_{\tilde y\rightarrow 0}\tilde y^{1-2\gamma}\partial_{\tilde y}\tilde u=(\frac{y}{\tilde y})^{n-s+2\gamma}v \;\;\text{ on } \;\;M.
\end{cases}
\end{equation*}
We may rephrase this via \;$\varphi=(\frac{\tilde y}{y})^{\frac{n-2\gamma}{2}}=(\frac{\tilde y}{y})^{n-s}$\; as
\begin{equation*}
\begin{cases}
D_{g}(\varphi u)=\varphi^{\frac{s+1}{n-s}} f\;\; \text{ in } \;\;X \\
\varphi u=\varphi v \;\,\;\text{ on } \;\,\;M
\end{cases}
\Longleftrightarrow
\begin{cases}
D_{\tilde g} u=f \;\,\;\text{in} \;\;X \\
u=v\;\;\text{on} \;\;M
\end{cases}
\end{equation*} 
and
\begin{equation*}
\begin{cases}
D_{g}(\varphi u)=\varphi^{\frac{s+1}{n-s}}f \;\;\text{in}\;\; X \\
\lim_{y\rightarrow 0}y^{1-2\gamma}\partial_{y}(\varphi u)=\varphi^{\frac{n+2\gamma}{n-2\gamma}}v \;\;\text{ on } \;\;M^{n}
\end{cases}
\Longleftrightarrow
\begin{cases}
D_{\tilde g}\tilde u=f \;\;\text{ in } \;\;X \\
\lim_{\tilde y\rightarrow 0}\tilde y^{1-2\gamma}\partial_{\tilde y} u=v \;\;\text{ on } \;\;M.
\end{cases}
\end{equation*}
Noticing \;$\frac{s+1}{n-s}=\frac{n+2+2\gamma}{n-2\gamma}$\; we thus have 
shown

\begin{equation*}
\xymatrix@R=6ex@C=-11ex
{
*+[l]
{
P^{\gamma}_{f,\tilde h}(\underline v)=\overline{v}
}
\ar@{<=>}[d] \ar@{<=>}[r]
&
*+[r]
{
{
\begin{cases}
D_{\tilde g}u=f \;\;\text{ in } \;\;X \\
u=\underline v \;\;\text{ on } \;\;M \\
-d_\gamma^*\lim_{\tilde y\rightarrow 0}\tilde y^{1-2\gamma}\partial_{\tilde y} u=\overline v \;\;\text{ on } \;\;M
\end{cases}
}
}
\ar@{<=>}[d]
\\
*+[l]
{
P^{\gamma}_{\varphi^{\frac{s+1}{n-s}}f,h }(\varphi \underline v)=\varphi^{\frac{n+2\gamma}{n-2\gamma}}\overline v
}
\ar@{<=>}[r]
&
*+[r]
{
{
\begin{cases}
D_{ g}(\varphi u)=\varphi^{\frac{n+2+2\gamma}{n-2\gamma}}f \text{ in } X \\
\varphi u=\varphi \underline v \text{ on } M \\
-d_\gamma^*\lim_{y\rightarrow 0}y^{1-2\gamma}\partial_{ y} (\varphi 
u)=\varphi^{\frac{n+2\gamma}{n-2\gamma}}\overline v \text{ on } M
\end{cases}
}
}
}
\end{equation*}
Therefore the non-homogeneous fractional operator verifies the conformal  property
$$
P^{\gamma}_{f,\tilde h}(\underline v)=\varphi^{-\frac{n+2\gamma}{n-2\gamma}}P^{\gamma}_{\varphi^{\frac{s+1}{n-s}}f,h }(\varphi \underline v)
\;\text{ for} \;
\tilde h=\varphi^{\frac{4}{n-2\gamma}}h
$$
or equivalently
$$P^{\gamma}_{\tilde f,\tilde h}(\underline v)
=
\varphi^{-\frac{n+2\gamma}{n-2\gamma}}P^{\gamma}_{f,h }(\varphi \underline v)
\; \text{ for } \;
\tilde h=\varphi^{\frac{4}{n-2\gamma}}h\;\text{ and }\;\tilde f=\varphi^{\frac{-s-1}{n-s}}f,
$$
hence extending the conformal property of  the homogeneous fractional operator to the non-homogeneous setting.  We remark that $$P^{\gamma}_{h}=P^{\gamma}_ {0,h}.$$

\section{Fundamental solutions in the asymptotically hyperbolic case}\label{fundso}
In this section, keeping the notations of the previous one, for an asymptotically hyperbolic manifold\;$(X, g^+)$ with conformal infinity 
\;$(M, [h])$, we study the existence and asymptotic behavior of the Poisson kernel \;$K_g:=K_{g}^{\gamma}$ \; of \;$D_g$, the Green's functions \;$\Gamma_g:=\Gamma_{g}^{\gamma}$\; of \;$D_g$\; under weighted normal boundary condition and \;$G_{h}:=G^{\gamma}_{h}$\; of the fractional conformal  Laplacian \;$P^{\gamma}_{h}$, i.e. 
\begin{equation*}\begin{split}
\begin{cases}
D_{g}K_g(\cdot, \xi)=0\;\; \text{ in }\;\; X \;\;\text{ and for all }\;\; \xi\in M\\
\lim_{y\rightarrow 0}K_g(y, x, \xi)=\d_{x}(\xi)\;\;\text{ and for all }\;\; x,\;\xi\in M
\end{cases}
\end{split}\end{equation*}
and
\begin{equation*}\begin{split}
\begin{cases}
D_{g}\Gamma_g(\cdot, \xi)=0\; \text{ in } X \;\;\text{ and for all }\;\; \xi\in M\\
-d_\gamma^*\lim_{y\rightarrow 0}y^{1-2\gamma}\partial_y\Gamma_g(y, x, \xi)=\d_{x}(\xi)\;\;\text{ and for all }\;\; x,\xi\in M,
\end{cases}
\end{split}\end{equation*}
where \;$d_\gamma^*$ is given by \eqref{dsgamma}, and 
$P^\gamma_h G^{\gamma}_{h}(x, \xi)=\d_{x}(\xi),\;\;x\in M.
$
So by definition
\begin{equation*}\begin{split}
K_{g}:(\overline X \times M)\setminus Diag(M)\longrightarrow\R_{+}
\end{split}\end{equation*}
is the Green's  function to the extension problem
\begin{equation*}\begin{split}
\begin{cases}
D_{g}U=0\; \text{ in } \;\;X \\
U=\underline{v}\;\; \text{on} \;\;M,
\end{cases}
\end{split}\end{equation*}
while
\begin{equation*}\begin{split}
\Gamma_{g}:(\overline X\times M)\setminus Diag(M)\longrightarrow\R
\end{split}\end{equation*}
is the Green's function to the dual problem
\begin{equation*}\begin{split}
\begin{cases}
D_{g}U=0\;\;\text{ in }\;\;X\\
-d_\gamma^*\lim_{y\to 0}y^{1-2\gamma}\partial_{y}U=\overline{v}  \;\;\text{ on } \;\;M
\end{cases}
\end{split}\end{equation*}
and 
\begin{equation*}\begin{split}
G_{h}:(M\times M)\setminus Diag(M)\longrightarrow\R.
\end{split}\end{equation*}
is the Green's function of the nonlocal problem
$
P^{\gamma}_{h} \underline{v}=\overline{v}  \;\;\text{ on } \;\;M.
$
They are linked via
\begin{equation}\label{regreen}
\Gamma_{g}=K_{g}*G_{h},
\end{equation}
where \;$*$\; denotes the standard convolution operation.

\subsection{Study of the Poisson kernel for \;$D_g$\; }
In this subsection we study the Poisson kernel \;$K_g$\; focusing on the existence issue and its asymptotics. We follow the method of Lee-Parker\cite{lp} of killing deficits successively. However,  due to the rigidity property involved in the problem, see  the normal form \eqref{eq:uniqdef},  we have to work close to the boundary in Fermi coordinates rather than normal ones. To compensate this we are forced to pass from the space of polynomials used in \cite{lp} to the space of homogeneous functions. We start with recalling  some related facts in the  case of the standard Euclidean space \;$\R^{n+1}_+$. According to \cite{cafsyl} on \;$\R^{n+1}_+$
\begin{equation}\label{Poisson_kernel_flat}
\begin{split}
K(y,x,\xi)=K^\gamma(y, x, \xi)= p_{n, \gamma}\frac{y^{2\gamma}}{(y^{2}+\vert x-\xi\vert^{2})^{\frac{n+2\gamma}{2}}},
\end{split}
\end{equation} 
where \;$p_{n, \gamma}$\; is as in \eqref{pngamma},  is the Poisson kernel of the operator
\begin{equation*}\label{Extension_Operator_flat}
D=-div(y^{1-2\gamma}\nabla (\,\cdot\,)),
\end{equation*}
namely the Green's function of the extension problem 
\begin{equation*}
\begin{cases}
Du=0 \;\;\text{ in }\;\; \R^{n+1}_+ \\
u=f \;\;\text{ on }\;\;\R^{n},
\end{cases}
\end{equation*} 
 i.e.
\begin{equation}\label{poisson_defining_equation}
\begin{cases}
DK(y,x,\xi)=0 \;\;\text{ in }\;\;\R^{n+1}_+\;\; \\
K(y,x,\xi)\rightarrow\delta_{x}(\xi)\; \;\;\text{ for }\;\; y\rightarrow  0.
\end{cases}
\end{equation} 
We will construct the Poisson kernel for \;$D_{g}$, cf.  \eqref{eqdg}, namely the Green's function of the analogous extension problem
\begin{equation*}
\begin{cases}
D_{g}u=0 \;\;\text{in} \;\;X \\
u=f \;\;\text{ on } \;\;M,
\end{cases}
\end{equation*} 
 i.e. \;$K_g$\; solves for \;$z\in X$\; and \;$\xi \in M$
\begin{equation*}
\begin{cases}
D_{g}K_{g}(z,\xi)=0\;\;\text{ in } \;\;X\\
K(z,\xi)\rightarrow\delta_{x}(\xi)\;\;\text{ for } \;\;y\rightarrow 0,
\end{cases}
\end{equation*}
where \;$z=(y,x)\in X$\; for \;$z$\; close to \;$M$. To that end we identify
\begin{equation*}
\xi \in M\cap U\subset U\cap X \text{ with }0\in B^{\R^{n+1}}_{\epsilon}(0)\cap 
\R^{n}\subset  B_{\epsilon}^{\R^{n+1}}(0)\cap \R^{n+1}_{+}
\end{equation*}
for some open neighborhood \;$U$\; of \;$\xi$\; in \;$X$\; and small \;$\epsilon>0$, and write \;$K(z)=K(z,0)$. We then have 
\begin{equation}\label{poisson_local_equation_for_the_poisson_kernel}
D_{g}K=-\frac{\partial_{p}}{\sqrt{g}}(\sqrt{g}g^{p,q}y^{1-2\gamma}\partial_{q}
K)+E_{g}K=f\in
y H_{-n-2\gamma-1}C^{\infty}
\end{equation}
on $\;B_{\epsilon}^{\R^{n+1}}(0)\cap \R^{n+1}_{+}$\; due \eqref{sqrt_g_sim_yC_infty}, which relies on minimality \;$H_{g}=0$,
where by definition 
\begin{equation}\label{definition_H_l}
H_{l}=\{\varphi\in C^{\infty}(\overline{\R^{n+1}_{+}}\setminus \{0\})\mid 
\varphi \text{ is homogeneous of degree } l\}.
\end{equation}
The next lemma allows us to solve homogeneous deficits 
homogeneously.
\begin{lem}
\label{poisson_homogeneous_solvability_dirichlet}$_{}$\\
For $\frac{1}{2}\neq \gamma \in (0, 1)$ and \;$f_{l}\in yH_{l-1}, \;l\in \N-n-2\gamma$\; there exists \;$K_{1+2\gamma+l}\in 
y^{2\gamma}H_{l+1}$\; such, that 
$$DK_{1+2\gamma+l}=f_{l}.$$
\end{lem}
\begin{pf}
First of all the Stone-Weierstra{\ss} Theorem implies
\begin{equation*}
\langle Q^{k}_{l}(y,x)=y^{2\gamma+2k}P_{l}(x)\mid k,l\in \N\text{ and }P_{l}\in 
\Pi_{l}\rangle \underset{\text{dense}}{\subset} y^{2\gamma}C^{0}(\overline 
B_{1}(0)\cap\R^{n+1}_{+})
\end{equation*}
and an easy induction argument shows, that we have a unique representation
\begin{equation*}
Q^{k}_{l}=\sum \vert z \vert^{2i}A_{2k+l-2i}
\end{equation*}
with \;$D$-harmonics of the form \;$A_{m}(y,x)=\sum y^{2\gamma+2l}P_{m-2l}(x),\;
DA_{m}=0$. Since 
\begin{equation*}
y^{2\gamma}C^{0}(S^{n}_{+})\underset{\text{dense}}{\subset} 
L^{2}_{y^{1-2\gamma}}(S^{n}_{+}),
\end{equation*}
we thus obtain a \;$D$-harmonic basis \;$E=\{e^{i}_{k}\}$\; for 
\;$L^{2}_{y^{1-2\gamma}}(S^{n}_{+})$\; with
\begin{equation*}
De^{i}_{k}=0, \;\;k=deg(e^{i}_{k}), \;\; k\in \N +2\gamma\text{ and } \;\;i\in 
\{1,\dots,d_{k}\},
\end{equation*}
where \;$d_{k}$\; denotes the dimension of the space of \;$D$-harmonics of degree 
$k$. We may assume, that \;$e^{i}_{k},e^{j}_{k}$ for $i\neq j$\; are orthogonal 
with respect to the scalar product on \;$L^{2}_{y^{1-2\gamma}}(S^{n}_{+})$. 
Moreover on \;$S^{n}_{+}$\; we have
\begin{equation*}
\begin{split}
0
= &
-De^{i}_{k}
=
\partial_{y}(y^{1-2\gamma}\partial_{y}e^{i}_{k})+y^{1-2\gamma}\Delta_{x}e^{i}_{k
}
=
\nabla y^{1-2\gamma}\nabla e^{i}_{k}+y^{1-2\gamma}\Delta e^{i}_{k} \\
= &
\nabla y^{1-2\gamma} \nabla 
e^{i}_{k}+y^{1-2\gamma}\frac{\Delta_{S^{n}}}{r^{2}}e^{i}_{k}+y^{1-2\gamma}[
\partial_{r}^{2}+\frac{n\partial_{r}}{r}]e^{i}_{k} \\
= &
\nabla_{S^{n}_{+}}^{\perp} y^{1-2\gamma}\nabla_{S^{n}_{+}}^{\perp} e^{i}_{k}
+
div_{S^{n}_{+}}(y^{1-2\gamma}\nabla_{S^{n}_{+}}e^{i}_{k})
+
k(k+n-1)y^{1-2\gamma}e^{i}_{k},
\end{split}
\end{equation*}
whence due to 
\begin{equation*}
\begin{split}
\nabla_{S^{n}_{+}}^{\perp} y^{1-2\gamma}\nabla_{S^{n}_{+}}^{\perp} e^{i}_{k}
= &
\langle \nabla y^{1-2\gamma}, \nu_{S^{n}_{+}}\rangle \langle 
\nu_{S^{n}_{+}}, \nabla e^{i}_{k}\rangle
=
(1-2\gamma)y^{-2\gamma}\langle e_{n+1}, \nu_{S^{n}_{+}}\rangle 
r\partial_{r}e^{i}_{k} 
= 
(1-2\gamma)ky^{1-2\gamma}e^{i}_{k}
\end{split}
\end{equation*}
there holds for\; 
$D_{S^{n}_{+}}=-div_{S^{n}_{+}}(y^{1-2\gamma}\nabla_{S^{n}_{+}}\,\cdot\,)$
\begin{equation*}
D_{S^{n}_{+}}e^{i}_{k}=k(k+n-2\gamma)y^{1-2\gamma}e^{i}_{k}.
\end{equation*}
Therefore \;$E=\{e^{i}_{k}\}$\; is an orthogonal basis of 
\;$y^{2\gamma-1}D_{S^{n}_{+}}$-eigenfunctions with eigenvalues 
$$\lambda_{k}=k(k+n-2\gamma).$$ 
By the same argument solving
\begin{equation}\label{poisson_homogeneous_to_solve}
\begin{cases}
Du=f\in L^{2}_{y^{2\gamma-1}}(\R^{n+1}_{+})\;\; \text{ in } \;\;\R^{n+1}_{+}\\
u=0  \;\;\text{ on } \;\; \R^{n}
\end{cases}
\end{equation}
with homogeneous \;$f, \;u$\; of degree \;$\lambda, \;\lambda+1+2\gamma$\; is equivalent to 
solving
\begin{equation*}
\begin{cases}
D_{S^{n}_{+}}u=f+(\lambda+1+2\gamma)(\lambda+n+1)y^{1-2\gamma}u\;\;\text{ in 
}\;\;S^{n}_{+} \\
u=0  \;\;\text{ on }\;\; \partial S^{n}_{+}=S^{n-1}
\end{cases}
\end{equation*}
and thus, writing 
$
u=\sum a_{i,k}e^{i}_{k}, \,y^{2\gamma-1}f=\sum 
b_{j,l}e^{j}_{l}
$
, also equivalent to solving
\begin{equation*}
\sum a_{i,k}(k(k+n-2\gamma)-(\lambda+1+2\gamma)(\lambda+n+1))e^{i}_{k}=\sum 
b_{j,l}e^{j}_{l}
\end{equation*}
and the latter system is always solvable in case
\begin{equation}\label{solvability_homogeneous_poisson}
k(k+n-2\gamma)-(\lambda+1+2\gamma)(\lambda+n+1)\neq 0 \;\;\text{ for all } \;\;
k, \;n, \;\lambda\in \N.
\end{equation}
This observation allows us to prove the lemma, by whose assumptions  
\begin{equation*}
deg(f_{l})=\lambda=m-n-2\gamma, \;\; m\in \N.
\end{equation*}
And we know 
\begin{equation*}
deg(e^{i}_{k})=k=m^{\prime}+2\gamma, \;\; m^{\prime}\in \N. 
\end{equation*}
Plugging these values into \eqref{solvability_homogeneous_poisson}, solvability 
of \eqref{poisson_homogeneous_to_solve} is a consequence of
\begin{equation*}
(m^{\prime}+2\gamma)(m^{\prime}+n)-(m-n+1)(m+1-2\gamma)\neq 0\;\;\text{ for all } \;\;
m^{\prime}, \; \;n,m\in \N
\end{equation*}
and this holds true for \;$\frac{1}{2}\neq\gamma \in (0,1)$. Thus we have 
proven solvability of
\begin{equation*}
\begin{cases}
DK_{1+2\gamma+l}=f_{l} \;\;\text{ in }\;\;\R^{n+1}_{+} \\
K_{1+2\gamma+l}=0  \;\;\text{ on } \;\;\R^{n}\setminus \{0\}
\end{cases}
\end{equation*}
with \;$K_{1+2\gamma+l}$\; being homogeneous of degree \;$1+2\gamma+l$. We are left 
with showing \;$K_{1+2\gamma+l}\in y^{2\gamma}H_{l+1}$. But this follows easily from 
Proposition \ref{poisson_regularity} below.
\end{pf}
\vspace{6pt}

\noindent
Now recalling \eqref{poisson_local_equation_for_the_poisson_kernel} we may use 
Lemma \ref{poisson_homogeneous_solvability_dirichlet} to solve 
\eqref{poisson_defining_equation} successively, since
\begin{equation*}
D_{g}K_{1+2\gamma+l}=f_{l}+(D_{g}-D)K_{1+2\gamma+l}\in f_{l}+ yH_{l}C^{\infty}
\end{equation*}
due to \eqref{sqrt_g_sim_yC_infty} and  \;$K_{1+2\gamma+l}\in y^{2\gamma}H_{l+1}$. With a suitable cut-off function
\begin{equation}\label{etaxi}
\eta_{\xi}:\overline X\longrightarrow \R^{+}, \;\;supp(\eta_{\xi})=B_{\epsilon}^{+}(\xi)=B_{\epsilon}^{g, +}(\xi)\;\;\text{ 
for } \;\;M\ni \xi \sim 0 \in \R^{n} \;\;\text{and}\;\;\epsilon>0\;\;\text{ small}
\end{equation}
and for the meaning of \;$B_{\epsilon}^{g, +}(\xi)$\; see Section \ref{notation_and_prelimiaries}, 
we then find   
\begin{equation*}
K_{g}=\eta_{\xi}(K+\sum_{l=-n-2\gamma}^{m+2-2\gamma} K_{1+2\gamma+l})+\kappa_{m} 
\end{equation*}
for \;$m\in \N$\; and a weak solution
\begin{equation*}
\begin{cases}
D_{g}\kappa_{m}=-D_{g}\left(\eta_{\xi}(K+\sum_{l=-n-2\gamma}^{m+2-2\gamma} 
K_{1+2\gamma+l})\right)=h_{m} \;\;\text{ in }\;\,\;X \\
\kappa_{m}=0\;\;\text{ on } \;\;M 
\end{cases}
\end{equation*}
with \;$h_{m} \in yC^{m,\alpha}$. 
The following weak regularity statement will be sufficient for our purpose.
\begin{pro}
\label{poisson_regularity} $_{}$\\
Let \;$h\in yC^{2k+3,\alpha}(X)$\; and \;$u\in W^{1,2}_{y^{1-2\gamma}}(X)$\;  be a 
weak solution of
\begin{equation*}
\begin{cases}
D_{g}u=h\;\;\text{ in }\;\;X \\
u=0\;\; \text{ on } \;\;M.
\end{cases}
\end{equation*}
Then \;$u$\; is of class \;$y^{2\gamma}C^{2k,\beta}(X)$, provided \;$H_{g}=0$.
\end{pro}
\vspace{6pt}

\noindent
Putting these facts together before giving the proof of Proposition \ref{poisson_regularity},  we have the existence of \;$K_g$\; and can describe its asymptotic.
\begin{cor}
\label{Poisson_kernel_asymptotics}$_{}$\\
Let \;$\frac{1}{2}\neq \gamma \in (0,1)$. Then  \;$K_g$\; exists and we may expand in \;$g$-normal 
Fermi-coordinates around $\xi\in M$
\begin{equation*}
K_{g}(z,\xi)
\;\in
\eta_{\xi}(z)
\left( p_{n, \gamma}
\frac{y^{2\gamma}}{\vert z 
\vert^{n+2\gamma}}+\sum^{2m+5-2\gamma}_{l=-n-2\gamma}y^{2\gamma}H_{1+l}(z)
\right)
+
y^{2\gamma}C^{2m,\alpha}(X)
\end{equation*}
with \;$H_{l}\in C^{\infty}(\R^{n+1}_+\setminus \{0\})$\; being homogeneous of order 
$l$\; and \;$p_{n, \gamma}$\; is as in \eqref{pngamma}, provided \;$H_{g}=0$.
\end{cor}
\vspace{6pt}

\noindent
\begin{pfn}{ of Proposition \ref{poisson_regularity}}
\vspace{2pt}

\noindent
We use the Moser iteration argument. First let 
$
p,q=1,\ldots,n+1\text{ and } i,j=1,\ldots,n
$
such, that $g_{n+1,i}=g_{y,i}=0$. The statement clearly holds by standard local 
regularity away from the boundary, since $D_{g}$ is strongly elliptic there. 
Now fixing a point $\xi\in M$ and a cut-off function
\begin{equation*}
\eta\in C^{\infty}_{0}(B^{+}_{r_{2}}(0),\R_{+}),\,\; \eta\equiv 1\text{ on } 
B^{+}_{r_{1}}(0)\;\;\text {for }\; 0<r_{1}<r_{2}\ll 1,\;\text{ where } \;\;M\ni \xi \sim 0\in \R^{n},
\end{equation*}
we pass to $g$-normal Fermi-coordinates around $\xi$ and estimate for some \;$\lambda\geq 2$\; and \;$\alpha\in \N^{n}$ 
\begin{equation}\label{poisson_kernel_moser_iteration_1}
\begin{split}
\underset{\R^{n+1}_{+}}{\int}y^{1-2\gamma}\vert \nabla_{z}(\vert 
\partial^{\alpha}_{x}u\vert^{\frac{\lambda}{2}}\eta)\vert^{2} 
\leq &
2\underset{\R^{n+1}_{+}}{\int}y^{1-2\gamma}\vert \nabla_{z}\vert 
\partial^{\alpha}_{x}u\vert^{\frac{\lambda}{2}}\vert^{2}\eta^{2}
+
2\underset{\R^{n+1}_{+}}{\int}y^{1-2\gamma}\vert \partial^{\alpha}_{x}u\vert^{\lambda}\vert 
\nabla_{z}\eta\vert^{2}
\end{split}
\end{equation}
and
\begin{equation*}
\begin{split}
\underset{\R^{n+1}_{+}}{\int}& y^{1-2\gamma}\vert \nabla_{z}\vert 
\partial^{\alpha}_{x}u\vert^{\frac{\lambda}{2}}\vert^{2}\eta^{2}
=
\frac{\lambda^{2}}{4}\underset{\R^{n+1}_{+}}{\int}y^{1-2\gamma} \nabla_{z} \partial^{\alpha}_{x}u
\nabla_{z} \partial^{\alpha}_{x}u \vert  
\partial^{\alpha}_{x}u\vert^{\lambda-2}\eta^{2}\\
= &
\frac{\lambda^{2}}{4(\lambda-1)}\underset{\R^{n+1}_{+}}{\int}y^{1-2\gamma} \nabla_{z} 
\partial^{\alpha}_{x}u
\nabla_{z} (\partial^{\alpha}_{x}u \vert  
\partial^{\alpha}_{x}u\vert^{\lambda-2}\eta^{2}) 
-
\frac{\lambda^{2}}{2(\lambda-1)}\underset{\R^{n+1}_{+}}{\int}y^{1-2\gamma} \nabla_{z} 
\partial^{\alpha}_{x}u
\partial^{\alpha}_{x}u \vert  
\partial^{\alpha}_{x}u\vert^{\lambda-2}\nabla_{z}\eta \eta \\
\leq &
\frac{\lambda^{2}}{4(\lambda-1)}\underset{\R^{n+1}_{+}}{\int}y^{1-2\gamma} \nabla_{z} 
\partial^{\alpha}_{x}u
\nabla_{z} (\partial^{\alpha}_{x}u \vert  
\partial^{\alpha}_{x}u\vert^{\lambda-2}\eta^{2}) \\
& +
\frac{\lambda^{2}}{8}\underset{\R^{n+1}_{+}}{\int}y^{1-2\gamma} \vert \nabla_{z} 
\partial^{\alpha}_{x}u\vert^{2}
\vert  \partial^{\alpha}_{x}u\vert^{\lambda-2}\eta^{2} 
+
\frac{\lambda^{2}}{2(\lambda-1)^{2}}\underset{\R^{n+1}_{+}}{\int}y^{1-2\gamma} 
 \vert  \partial^{\alpha}_{x}u\vert^{\lambda}\vert \nabla_{z}\eta \vert^{2}.
\end{split}
\end{equation*}
Absorbing the second summand above this implies
\begin{equation}\label{poisson_kernel_moser_iteration_3}
\begin{split}
\underset{\R^{n+1}_{+}}{\int}
\hspace{-2pt}
y^{1-2\gamma}&\vert \nabla_{z}(\vert 
\partial^{\alpha}_{x}u\vert^{\frac{\lambda}{2}})\vert^{2}\eta^{2} 
\leq 
\frac{\lambda^{2}}{2(\lambda-1)}\underset{\R^{n+1}_{+}}{\int}
\hspace{-2pt}
D(\partial^{\alpha}_{x}u)
\partial^{\alpha}_{x}u \vert  \partial^{\alpha}_{x}u\vert^{\lambda-2}\eta^{2} 
+
\frac{\lambda^{2}}{(\lambda-1)^{2}}\underset{\R^{n+1}_{+}}{\int}
\hspace{-2pt}
y^{1-2\gamma}\vert 
\partial^{\alpha}_{x}u\vert^{\lambda}\vert \nabla_{z}\eta\vert^{2}
\end{split}
\end{equation}
Due to
$
D(\partial^{\alpha}_{x}u)=\partial^{\alpha}_{x}Du,
$
and the structure of the metric 
\begin{equation}\label{poisson_kernel_moser_iteration_4}
\begin{split} 
\underset{\R^{n+1}_{+}}{\int} &\partial_{x}^{\alpha}(D  u)\partial^{\alpha}_{x}u\vert 
\partial^{\alpha}_{x}u\vert^{\lambda-2}\eta^{2} 
= 
\underset{\R^{n+1}_{+}}{\int} \partial_{x}^{\alpha}(D_{g}u)\partial^{\alpha}_{x}u\vert 
\partial^{\alpha}_{x}u\vert^{\lambda-2}\eta^{2} 
-
\underset{\R^{n+1}_{+}}{\int} \partial_{x}^{\alpha}((D_{g}-D)u)\partial^{\alpha}_{x}u\vert 
\partial^{\alpha}_{x}u\vert^{\lambda-2} \eta^{2} \\
= &
\underset{\R^{n+1}_{+}}{\int}
\partial_{x}^{\alpha}
[
h
+
\frac{\partial_{p}\sqrt{g}}{\sqrt{g}}y^{1-2\gamma}g^{p,q}\partial_{q}u
+
y^{1-2\gamma}\partial_{i}((g^{i,j}-\delta^{i,j})\partial_{j}u) \\
& \quad \quad\quad\quad\quad\quad\quad\quad\quad -
\frac{n-2\gamma}{2}\frac{\partial_{y}\sqrt{g}}{\sqrt{g}}y^{-2\gamma}u
]
\partial^{\alpha}_{x}u\vert \partial^{\alpha}_{x}u\vert^{\lambda-2}\eta^{2} 
= 
I_{1}+\ldots+I_{4}.
\end{split}
\end{equation}
We may assume \;$\nabla^{k}_{z}\eta\leq \frac{C}{\epsilon^{k}}$\; for \;$k=0,1,2$, 
where $\epsilon=r_{2}-r_{1}$. Then
\begin{enumerate}[label=(\roman*)]
\item \quad 
$
\vert I_{1} \vert
\leq  
C\int_{\R^{n+1}_{+}}\vert \nabla ^{\vert \alpha \vert}_{x}h\vert  
\vert \nabla^{\vert \alpha \vert}_{x}u\vert^{\lambda-1} \eta^{2}
$

\item using integrations by parts and \eqref{sqrt_g_sim_yC_infty}
\begin{equation*}
\begin{split}
\vert I_{2}\vert
\leq &
\vert 
\underset{\R^{n+1}_{+}}{\int}
y^{1-2\gamma}\partial_{q}\partial^{\alpha}_{x}(\frac{\partial_{p}\sqrt{g}}{\sqrt
{g}}g^{p,q}u)\partial^{\alpha}_{x}u\vert \partial^{\alpha}_{x}u\vert 
^{\lambda-2} \eta^{2}\vert +
\vert \underset{\R^{n+1}_{+}}{\int} 
y^{1-2\gamma}\partial^{\alpha}_{x}(\partial_{q}(\frac{\partial_{p}\sqrt{g}}{
\sqrt{g}}g^{p,q})u)\partial^{\alpha}_{x}u\vert \partial^{\alpha}_{x}u\vert 
^{\lambda-2} \eta^{2}\vert \\
\leq &
\vert \underset{\R^{n+1}_{+}}{\int} 
y^{1-2\gamma}\partial_{y}\partial^{\alpha}_{x}(\frac{\partial_{y}\sqrt{g}}{\sqrt
{g}}u)\partial^{\alpha}_{x}u\vert \partial^{\alpha}_{x}u\vert ^{\lambda-2} 
\eta^{2} \vert   
+
\frac{C_{\vert \alpha \vert}}{\lambda}\sum_{m\leq \vert \alpha \vert}\underset{\R^{n+1}_{+}}{\int}
y^{1-2\gamma}\vert \nabla_{x}^{m}u\vert \vert 
\partial^{\alpha}_{x}u\vert^{\frac{\lambda-2}{2}}
\vert \nabla_{x} \vert \partial^{\alpha}_{x}u\vert^{\frac{\lambda}{2}}\vert 
\eta^{2} \\
& +
C_{\vert \alpha \vert}\sum_{m\leq \vert \alpha \vert} \underset{\R^{n+1}_{+}}{\int} y^{1-2\gamma}\vert 
\nabla_{x}^{m}u\vert \vert \partial^{\alpha}_{x}u\vert^{\lambda-1} [\vert 
\nabla_{x}\eta\vert \eta+\eta^{2}] \\
\leq &
\frac{C_{\vert \alpha \vert}}{\lambda}\sum_{m\leq \vert \alpha \vert} 
\underset{\R^{n+1}_{+}}{\int}y^{1-2\gamma}\vert \nabla_{x}^{m}u\vert^{\frac{\lambda}{2}} 
\vert \nabla_{z} \vert \partial^{\alpha}_{x}u\vert^{\frac{\lambda}{2}}\vert 
\eta^{2} 
+
C_{\vert \alpha \vert}\sum_{m\leq \vert \alpha \vert} \underset{\R^{n+1}_{+}}{\int}y^{1-2\gamma}\vert 
\nabla_{x}^{m}u \vert^{\lambda} [\vert \nabla_{z}\eta\vert \eta+\eta^{2}] \\
\end{split}
\end{equation*}

\item using integration by parts and recalling $i,j=1,\ldots,n$
\begin{equation*}
\begin{split}
\vert I_{3} \vert
\leq &
\frac{C}{\lambda} \underset{\R^{n+1}_{+}}{\int}y^{1-2\gamma}\vert 
\partial^{\alpha}_{x}((g^{i,j}-\delta^{i,j})\partial_{j}u)\vert \vert 
\partial^{\alpha}_{x}u\vert^{\frac{\lambda-2}{2}}\vert \partial_{i}\vert 
\partial^{\alpha}_{x}u\vert^{\frac{\lambda}{2}}\vert \eta^{2} \\
& +
C\underset{\R^{n+1}_{+}}{\int}y^{1-2\gamma}\vert 
\partial^{\alpha}_{x}((g^{i,j}-\delta^{i,j})\partial_{j}u)\vert \vert 
\partial^{\alpha}_{x}u\vert^{\lambda-1}\vert \partial_{i}\eta \vert \eta \\
\leq &
\frac{C}{\lambda^{2}}\sup_{B^{+}_{r_{2}}}\vert g^{i,j}-\delta^{i,j}\vert \underset{\R^{n+1}_{+}}{\int}
y^{1-2\gamma}
\vert \partial_{i}\vert \partial^{\alpha}_{x}u\vert^{\frac{\lambda}{2}}\vert
\vert \partial_{j}\vert \partial^{\alpha}_{x}u\vert^{\frac{\lambda}{2}}\vert
\eta^{2}  \\
& +
\frac{C_{\vert \alpha \vert}}{\lambda}\sum_{m\leq \vert \alpha }
 \underset{\R^{n+1}_{+}}{\int}y^{1-2\gamma}\vert \nabla_{x}^{m}u\vert \vert 
\partial^{\alpha}_{x}u\vert^{\frac{\lambda-2}{2}}\vert \nabla_{x}\vert 
\partial^{\alpha}_{x}u\vert^{\frac{\lambda}{2}}\vert \eta^{2}\\
& +
\frac{C}{\lambda}\sup_{B^{+}_{r_{2}}}\vert g^{i,j}-\delta^{i,j}\vert \underset{\R^{n+1}_{+}}{\int}
y^{1-2\gamma} \vert \partial_{j}\vert 
\partial^{\alpha}_{x}u\vert^{\frac{\lambda}{2}}\vert \vert 
\partial^{\alpha}_{x}u\vert^{\frac{\lambda}{2}}\vert\partial_{i} \eta \vert 
\eta 
+
C_{\vert \alpha \vert}\sum_{m\leq \vert \alpha }
 \underset{\R^{n+1}_{+}}{\int}y^{1-2\gamma}\vert \nabla_{x}^{m}u\vert^{\lambda} \vert 
\nabla_{x}\eta\vert \eta 
\end{split}
\end{equation*}

\item \quad 
$
\vert I_{4} \vert 
\leq 
C_{\vert \alpha \vert}\sum_{m\leq \vert\alpha\vert}\underset{\R^{n+1}_{+}}{\int}y^{1-2\gamma} \vert 
\nabla^{m}_{x}u\vert^{\lambda}
 \eta^{2}
$
using \eqref{sqrt_g_sim_yC_infty}.
\end{enumerate}
Applying H\"older's and Young's inequality to (i)-(vi)  we obtain
\begin{equation*}
\begin{split}
\sum^{4}_{i=1}\vert I_{i}\vert
\leq &
\frac{C\sup_{B^{+}_{r_{2}}}\vert g-\delta\vert}{\lambda^{2}}
\underset{\R^{n+1}_{+}}{\int}y^{1-2\gamma}\vert \nabla_{z}\vert \partial^{\alpha}_{x}u 
\vert^{\frac{\lambda}{2}}\vert^{2}\eta^{2} 
+
\frac{C_{\vert \alpha\vert}}{\epsilon^{2}}\sum_{k\leq \vert \alpha \vert}
\Vert 
\nabla^{k}_{x}u\Vert_{L^{\lambda}_{y^{1-2\gamma}}(B^{+}_{r_{2}})}^{\lambda} \\
& +
C\Vert y^{2\gamma-1}\nabla^{\vert \alpha 
\vert}_{x}h\Vert_{L^{\lambda}_{y^{1-2\gamma}}(B^{+}_{r_{2}})}
\Vert \nabla^{\vert \alpha 
\vert}_{x}u\Vert_{L^{\lambda}_{y^{1-2\gamma}}(B^{+}_{r_{2}})}^{\lambda-1}
.
\end{split}
\end{equation*}
We may assume \;$C\sup_{B^{+}_{r_{2}}}\vert g-\delta\vert<\frac{1}{2}$, whence in 
view of \eqref{poisson_kernel_moser_iteration_3} and 
\eqref{poisson_kernel_moser_iteration_4} 
\begin{equation*}
\begin{split}
\underset{\R^{n+1}_{+}}{\int}y^{1-2\gamma}\vert \nabla_{z}\vert 
\partial^{\alpha}_{x}u\vert^{\frac{\lambda}{2}}\vert^{2}\eta^{2}
\leq &
\frac{C_{\vert \alpha \vert}\lambda}{\epsilon^{2}}
\sum_{k\leq \vert \alpha \vert}
\Vert 
\nabla^{k}_{x}u\Vert_{L^{\lambda}_{y^{1-2\gamma}}(B^{+}_{r_{2}})}^{\lambda}
\\
& +
C\lambda\Vert y^{2\gamma-1}\nabla^{\vert \alpha 
\vert}_{x}h\Vert_{L^{\lambda}_{y^{1-2\gamma}}(B^{+}_{r_{2}})}
\Vert \nabla^{\vert \alpha 
\vert}_{x}u\Vert_{L^{\lambda}_{y^{1-2\gamma}}(B^{+}_{r_{2}})}^{\lambda-1},
\end{split}
\end{equation*} 
so \eqref{poisson_kernel_moser_iteration_1} implies
\begin{equation}\label{poisson_kernel_moser_iteration_2}
\begin{split}
\underset{\R^{n+1}_{+}}{\int}y^{1-2\gamma}\vert \nabla_{z}(\vert \partial^{\alpha}_{x} 
u\vert^{\frac{\lambda}{2}}\eta)\vert^{2}
\leq &
\frac{C_{\vert \alpha \vert}\lambda}{\epsilon^{2}}
\sum_{k\leq \vert \alpha \vert}
\Vert 
\nabla^{k}_{x}u\Vert_{L^{\lambda}_{y^{1-2\gamma}}(B^{+}_{r_{2}})}^{\lambda}
\\
& +
C\lambda\Vert y^{2\gamma-1}\nabla^{\vert \alpha 
\vert}_{x}h\Vert_{L^{\lambda}_{y^{1-2\gamma}}(B^{+}_{r_{2}})}
\Vert \nabla^{\vert \alpha 
\vert}_{x}u\Vert_{L^{\lambda}_{y^{1-2\gamma}}(B^{+}_{r_{2}})}^{\lambda-1}.
\end{split}
\end{equation} 
The weighted Sobolev inequality  of Fabes-Kenig-Seraponi \cite{fabkenser}  Theorem 1.2 with \;$\kappa=\frac{n+1}{n}$\; 
then shows
\begin{equation*}
\begin{split}
r_{2}^{-\frac{n+2\gamma}{n+1}}\Vert \partial^{ \alpha }_{x}u \Vert_{L^{\kappa 
\lambda}_{y^{1-2\gamma}}(B^{+}_{r_{1}})}^{\lambda}
\leq &
\frac{C_{\vert \alpha \vert}\lambda}{\epsilon^{2}}
\sum_{k\leq \vert \alpha \vert}
\Vert 
\nabla^{k}_{x}u\Vert_{L^{\lambda}_{y^{1-2\gamma}}(B^{+}_{r_{2}})}^{\lambda}
\\
& +
C\lambda\Vert y^{2\gamma-1}\nabla^{\vert \alpha 
\vert}_{x}h\Vert_{L^{\lambda}_{y^{1-2\gamma}}(B^{+}_{r_{2}})}
\Vert \nabla^{\vert \alpha 
\vert}_{x}u\Vert_{L^{\lambda}_{y^{1-2\gamma}}(B^{+}_{r_{2}})}^{\lambda-1}
.
\end{split}
\end{equation*}
By rescaling we may assume for some \;$0<\epsilon_{0}\ll 1$, that
\begin{equation}\label{poisson_rescaling_deficit}
\Vert u \Vert_{L^{2}_{y^{1-2\gamma}}}+\sum^{\vert \alpha \vert}_{k=0}\Vert  
y^{2\gamma-1}\nabla^{k}_{x}h\Vert_{L^{\infty}_{y^{1-2\gamma}}(B^{+}_{(2+\vert 
\alpha \vert)\epsilon_{0}})}= 1,
\end{equation}
and putting \;$\lambda_{i}=2(\frac{n+1}{n})^{i}$\; and 
\;$\rho_{i}=\epsilon_{0}(1+\frac{1}{2^{i}})$\; we obtain
\begin{equation*}
\begin{split}
\Vert \nabla^{ \vert \alpha \vert }_{x}u 
\Vert_{L^{\lambda_{i+1}}_{y^{1-2\gamma}}(B^{+}_{\rho_{i+1}})}
\leq 
\sqrt[\lambda_{i}]{C_{\vert\alpha\vert,\epsilon_{0}}\lambda_{i}2^{2i}}
\cdot
\sup_{m\leq \vert \alpha \vert}
[
 & \Vert \nabla^{m}_{x}u\Vert_{L^{\lambda_{i}}_{y^{1-2\gamma}}(B^{+}_{r_{2}})} 
+
\Vert 
\nabla^{m}_{x}u\Vert_{L^{\lambda_{i}}_{y^{1-2\gamma}}(B^{+}_{r_{2}})}^{\frac{1}{
2}}
],
\end{split}
\end{equation*}
where we have used \;$\frac{1}{2}\leq \frac{\lambda_{i}-1}{\lambda_{i}}<1$. 
Iterating this inequality then shows
\begin{equation*}
\begin{split}
\Vert \nabla^{ \vert \alpha  \vert}_{x}u 
\Vert_{L^{\infty}_{y^{1-2\gamma}}(B^{+}_{\epsilon_{0}})}
\leq 
C_{\alpha,\epsilon_{0}}(1+\sup_{m\leq \vert \alpha \vert}
\Vert \nabla^{m}_{x}u \Vert_{L^{2}_{y^{1-2\gamma}}(B^{+}_{2\epsilon_{0}})})
\leq 
C_{\alpha,\epsilon_{0}},
\end{split}
\end{equation*}
where the last inequality follows from iterating 
\eqref{poisson_kernel_moser_iteration_2} with \;$\lambda=2$
and \eqref{poisson_rescaling_deficit}. Rescaling back we conclude
\begin{equation}\label{poisson_kernel_moser_iteration_5}
\begin{split}
\sum^{m}_{k=0}\Vert \nabla^{ k }_{x}u 
\Vert_{L^{\infty}_{y^{1-2\gamma}}(B^{+}_{\epsilon_{0}})}
\leq 
C_{m,\epsilon_{0}}
[
\Vert u \Vert_{L^{2}_{y^{1-2\gamma}}}
+
\sum^{m}_{k=0}\Vert y^{2\gamma-1}\nabla^{ k }_{x}h 
\Vert_{L^{\infty}_{y^{1-2\gamma}}}
].
\end{split}
\end{equation}
Note, that
$
D(\partial^{\alpha}_{x}u)
= 
\partial_{x}^{\alpha}h
+
\partial^{\alpha}_{x}((D-D_{g})u),
$
where
\begin{equation}\label{poisson_x_deriviative_versus_difference_of_operators}
\begin{split}
\partial^{\alpha}_{x}((D& -D_{g})u)
= 
\partial^{\alpha}_{x}
[
\frac{\partial_{p}\sqrt{g}}{\sqrt{g}}y^{1-2\gamma}g^{p,q}\partial_{q}u 
+
\partial_{i}(y^{1-2\gamma}(g^{i,j}-\delta^{i,j})\partial_{j}u)
-
\frac{n-2\gamma}{2}\frac{\partial_{y}\sqrt{g}}{\sqrt{g}}y^{-2\gamma}u
]
\\
= &
\partial_{q}\partial_{x}^{\alpha}(\frac{\partial_{p}\sqrt{g}}{\sqrt{g}}y^{
1-2\gamma}g^{p,q}u)
-
\partial^{\alpha}_{x}(\partial_{q}(\frac{\partial_{p}\sqrt{g}}{\sqrt{g}}y^{
1-2\gamma}g^{p,q})u) \\
& +
\partial_{i}\partial^{\alpha}_{x}(y^{1-2\gamma}(g^{i,j}-\delta^{i,j})\partial_{j
}u)
-
\frac{n-2\gamma}{2}\partial^{\alpha}_{x}(\frac{\partial_{y}\sqrt{g}}{\sqrt{g}}y^
{-2\gamma}u).
\end{split}
\end{equation}
In particular, since 
$
-\partial_{p}(y^{1-2\gamma}g^{p,q}\partial_{q}v) 
=
Dv
-
\partial_{i}(y^{1-2\gamma}(g^{i,j}-\delta^{i,j})v)
$ 
we may write
\begin{equation*}
\partial_{p}(y^{1-2\gamma}g^{p,q}\partial_{q}\partial^{\alpha}_{x}u)
=
\partial^{\alpha}_{x}h
+
h^{\alpha}
+
\sum \partial_{p}h^{\alpha}_{p},
\end{equation*}
where \;$h^{\alpha}, \;h^{\alpha}_{p}$\; depend only on \;$x-$derivatives of \;$u$\; of 
order up to \;$\vert \alpha \vert$, and
due to \eqref{sqrt_g_sim_yC_infty}, \;\eqref{poisson_kernel_moser_iteration_5},  
there holds 
\begin{equation*}
\sum^{m}_{\vert \alpha \vert=0}\Vert 
\frac{h^{\alpha}}{y^{1-2\gamma}}, \frac{h^{\alpha}_{p}}{y^{1-2\gamma}}\Vert_{L^{
\infty}_{y^{1-2\gamma}}(B^{+}_{\epsilon_{0}})}
\leq
C_{m,\epsilon_{0}}
[
\Vert u \Vert_{L^{2}_{y^{1-2\gamma}}}
+
\sum^{m}_{k=0}\Vert y^{2\gamma-1}\nabla^{ k }_{x}h 
\Vert_{L^{\infty}_{y^{1-2\gamma}}}
]
\;\;\text{for all}\;\; m\in \N.
\end{equation*} 
Then Zamboni\cite{Zamboni} Theorem 5.2 shows H\"older regularity, i.e. for all \;$m\in \N$
\begin{equation}\label{poisson_zamboni}
\sum^{m}_{k=0}\Vert 
\nabla^{k}_{x}u\Vert_{C^{0,\alpha}(B^{+}_{\frac{\epsilon_{0}}{2}})}
\leq 
C_{m,\epsilon_{0}}
[
\Vert u \Vert_{L^{2}_{y^{1-2\gamma}}}
+
\sum^{m}_{k=0}\Vert y^{2\gamma-1}\nabla^{ k }_{x}h 
\Vert_{L^{\infty}_{y^{1-2\gamma}}}
].
\end{equation} 
This allows us to integrate the equation directly. Indeed from 
\eqref{poisson_x_deriviative_versus_difference_of_operators} we have
\begin{equation*}
D(\partial^{\alpha}_{x}u)
=
\partial^{\alpha}_{x}h
+
\partial^{\alpha}_{x}((D-D_{g})u)
=
\partial^{\alpha}_{x}h
+
\partial_{y}(y^{2-2\gamma}f^{\alpha}_{1})+y^{1-2\gamma}f^{\alpha}_{2},
\end{equation*}
where by definition 
\;$
f^{\alpha}_{1}=\partial^{\alpha}_{x}(\frac{\partial_{y}\sqrt{g}}{y\sqrt{g}}u)
$\;
and 
\begin{equation*}
\begin{split}
f_{2}^{\alpha}
= &
\partial_{i}\partial_{x}^{\alpha}(\frac{\partial_{j}\sqrt{g}}{\sqrt{g}}g^{i,j}u)
-
\partial^{\alpha}_{x}(y^{2\gamma-1}\partial_{q}(\frac{\partial_{p}\sqrt{g}}{
\sqrt{g}}y^{1-2\gamma}g^{p,q})u) \\
& +
\partial_{i}\partial^{\alpha}_{x}((g^{i,j}-\delta^{i,j})\partial_{j}u)
-
\frac{n-2\gamma}{2}\partial^{\alpha}_{x}(\frac{\partial_{y}\sqrt{g}}{y\sqrt{g}}
u)
.
\end{split}
\end{equation*}
This implies 
\begin{equation*}
-\partial_{y}(y^{1-2\gamma}\partial_{y}\partial^{\alpha}_{x}u)
=
\partial_{y}(y^{2-2\gamma}f^{\alpha}_{1})
+
y^{1-2\gamma}(f^{\alpha}_{2}+\Delta_{x}\partial^{\alpha}_{x}u+y^{2\gamma-1}
\partial^{x}_{\alpha}h)
\end{equation*}
and we obtain
\begin{equation}\label{poisson_equation_integrated_1}
\begin{split}
\partial^{\alpha}_{x}u(y,x)
= 
y^{2\gamma}\bar u_{0}^{\alpha}(x) 
- 
\int^{y}_{0}\sigma \tilde f^{\alpha}_{1}(\sigma,x)d\sigma
- 
\int^{y}_{0}\sigma^{2\gamma-1}\int^{\sigma}_{0}\tau^{1-2\gamma} \tilde 
f^{\alpha}_{2}(\tau,x)d\tau d\sigma,
\end{split}
\end{equation}
where by definition we may write with smooth coefficients \;$f_{i,\beta}$
\begin{equation}\label{poisson_tilde_f_structure}
\tilde f^{\alpha}_{1}=\sum _{\vert \beta\vert \leq \vert \alpha 
\vert}f_{1,\beta}\partial^{\beta}_{x}u
\;\text{ and }\;
\tilde f^{\alpha}_{2}=\frac{\partial^{\alpha}_{x}h}{y^{1-2\gamma}}+\sum _{\vert 
\beta\vert \leq \vert \alpha \vert+2}f_{2,\beta} \partial_{x}^{\beta}u.
\end{equation}
Let \;$h\in yC^{l,\lambda^{\prime}}$. Then \eqref{poisson_zamboni} shows
\begin{equation*}
\fa \vert \alpha \vert\leq l \;: \;\nabla^{\vert \alpha\vert}_{x}u \in 
C^{0,\lambda},
\end{equation*}
whence \;$\fa\; \vert \alpha\vert \leq l-2 \;: \;\tilde  
f_{i}^{\alpha}\in C^{0,\lambda}$\; due to \eqref{poisson_tilde_f_structure}. In 
particular  \eqref{poisson_equation_integrated_1} implies
\begin{equation*}
\partial^{\alpha}_{x}u(y,x)
=
y^{2\gamma}\bar u_{0}^{\alpha}(x) +o(y^{2\gamma}),
\end{equation*}
so \;$\bar u_{0}^{\alpha}\in C^{l+2,\lambda}$\; anyway by interior regularity. We define
\begin{equation}\label{poisson_y_minus_2_gamma_redefinition}
\bar u^{\alpha}=y^{-2\gamma}\partial^{\alpha}_{x}u
,\;\;
\bar f_{1}^{\alpha}=y^{-2\gamma}f_{1}^{\alpha}
,\;\;
\bar f_{2}^{\alpha}=y^{-2\gamma}\tilde f^{\alpha}_{2}.
\end{equation} 
We then find from  \eqref{poisson_equation_integrated_1}, that 
\begin{equation}\label{poisson_equation_integrated}
\begin{split}
\bar u^{\alpha}(y,x)
= &
\bar u_{0}^{\alpha}(x) 
-
y^{-2\gamma}\int^{y}_{0}\sigma^{1+2\gamma} \bar f^{\alpha}_{1}(\sigma,x) d\sigma
-
y^{-2\gamma}\int^{y}_{0}\sigma^{2\gamma-1}\int^{\sigma}_{0}\tau\bar 
f^{\alpha}_{2}(\tau,x)d\tau d\sigma\\
= &
\bar u_{0}^{\alpha}(x) + \bar u_{1}^{\alpha}(y,x) +\bar u_{2}^{\alpha}(y,x),
\end{split}
\end{equation}
where according to 
\eqref{poisson_tilde_f_structure}, \eqref{poisson_y_minus_2_gamma_redefinition} 
we may write with smooth coefficients \;$f_{i,\beta}$
\begin{equation}\label{poisson_bar_f_structure}
\bar f^{\alpha}_{1}=\sum _{\vert \beta\vert \leq \vert \alpha 
\vert}f_{1,\beta}\bar u^{\beta}
\;\;\text{ and }\;\;
\bar f^{\alpha}_{2}=\frac{\partial^{\alpha}_{x}h}{y}+\sum _{\vert \beta\vert 
\leq \vert \alpha \vert+2}f_{2,\beta}\bar u^{\beta}.
\end{equation}
Then \eqref{poisson_equation_integrated_1} and \;$\fa \;\vert 
\alpha\vert \leq l-2 \;: \;\tilde  f_{i}^{\alpha}\in C^{0,\lambda}$\;  already show
\begin{equation*}
\fa \vert \alpha\vert \leq l-2 \; :\;\bar u^{\alpha}\in C^{0,\lambda}
\end{equation*}
and we may assume 
$
\fa \vert \alpha\vert\leq  l-2-2m\;:\;\partial^{2m}_{y}\bar u^{\alpha}\in 
C^{0,\lambda}
$ 
 inductively, whence according to \eqref{poisson_bar_f_structure} 
\begin{equation*}
\fa \vert \alpha \vert\leq l-2-2(m+1) \;:\;\partial^{2m}_{y}\bar f^{\alpha}_{i}\in 
C^{0,\lambda}. 
\end{equation*}
Then \eqref{poisson_equation_integrated} implies via Taylor expansion
\begin{equation*}
\fa \vert \alpha \vert\leq l-2-2(m+1)\; :\;\partial_{y}^{2m+2}\bar 
u^{\alpha}_{i}, \;\partial_{y}^{2m+2}\bar u^{\alpha}\in C^{0,\alpha}.
\end{equation*}
Thus we have proven
$\fa\; \vert \alpha \vert\leq l-2-2m \; : \;\partial^{2m}_{y}\partial^{\alpha}_{x} 
u\in C^{0,\lambda}$\; for some \;$\lambda>0$. However, since there are only even 
powers in the \;$y$-derivative, we only find 
\;$u\in C^{l-3,\lambda}$\; for $\;l\in 2\N$.
The proof is thereby complete.\end{pfn}
\subsection{Green's function for \;$D_g$\; under weighted Neumann boundary condition}\label{subsec:Greens_function}
In this subsection we study the Green's function \;$\Gamma_g$. As in the previous one we consider the existence  and  asymptotics issue. To do that we use the method of  Lee-Parker\cite{lp} and have  the same difficulties to overcome as in the previous subsection. We first note that on \;$\R^{n+1}_+$
\begin{equation}\label{gngamma}
\Gamma(y,x,\xi)=\Gamma^\gamma(y,x,\xi)=\frac{g_{n, \gamma}}{(y^{2}+\vert x-\xi\vert^{2})^{\frac{n-2\gamma}{2}}}, \;\; (y, x)\in \R^{n+1}_+, \;\;\xi\in \R^n
\end{equation}
for some \;$g_{n, \gamma}>0$\; is the Green's function to the dual problem 
\begin{equation*}
\begin{cases}
Du=0\;\;\text{ in }\;\;\R^{n+1}_+ \\
-d_\gamma^*\lim_{y\to 0}y^{1-2\gamma}\partial_{y}u(y,\cdot)=f\;\; \text{ on }\;\; \R^{n},
\end{cases}
\end{equation*}
i.e.
\begin{equation*}
\begin{cases}
D\Gamma(, \xi)=0\;\;\text{ in }\;\;\R^{n+1}_+, \;\;\xi\in \R^n \\
-d_\gamma^*\lim_{y\to 0}y^{1-2\gamma}\Gamma(y,x,\xi)=\delta_{x}(\xi), \;\;x, \;\xi\in \R^n.
\end{cases}
\end{equation*}
We will construct the Green's function \;$\Gamma_g$\; for the analogous problem
\begin{equation*}
\begin{cases}
D_{g}u=0\;\; \text{ in }\;\;X \\
-d_\gamma^*\lim_{y\to 0}y^{1-2\gamma}\partial_{y}u(y,\cdot)=f \;\;\text{ on }\;\;M
\end{cases}
\end{equation*} 
for \;$D_{g}=-div_{g}(y^{1-2\gamma}\nabla_{g} (\,\cdot\,))+E_{g}$, i.e. for \;$z\in X$\; and \;$\xi \in M$
\begin{equation}\label{green_defining_equation} 
\begin{cases}
D_{g}\Gamma_{g}(\cdot,\xi)=0 \;\;\text{ in }\;\; X\\
-d_\gamma^*\lim_{y\to 0}y^{1-2\gamma}\Gamma_{g}(z,\xi)=\delta_{x}(\xi),
\end{cases}
\end{equation}
where \;$z=(y,x)\in X$\;in $g$-normal Fermi-coordinates close to \;$M$. To that end  we identify
\begin{equation*}
\xi \in M\cap U\subset U\cap X \;\;\text{ with }\;\;0\in B_{\epsilon}^{\R^{n+1}}(0)\cap 
\R^{n}\subset  B_{\epsilon}^{\R^{n+1}}(0)\cap \R^{n+1}_{+}
\end{equation*}
as in the previous subsection, and write \;$\Gamma(z)=\Gamma(z,0)$. On \;$B_{\epsilon}^{\R^{n+1}}(0)\cap \R^{n+1}_{+}$\; we 
then have 
\begin{equation}\label{green_local_equation_for_the_greens_function}
D_{g}\Gamma=-\frac{\partial_{p}}{\sqrt{g}}(\sqrt{g}g^{p,q}y^{1-2\gamma}\partial_
{ q}\Gamma)+E_{g}\Gamma
=
f
\in 
y^{1-2\gamma}H_{-n+2\gamma-1}C^{\infty}.
\end{equation}
Again we may solve homogeneous deficits homogeneously.
\begin{lem}
\label{green_homogeneous_solvability_neumann}$_{}$\\
For $\frac{1}{2}\neq \gamma \in (0,1)$ and $f_{l}\in y^{1-2\gamma}H_{l+2\gamma-1}, \; l\in \N-n$\; there exists 
\;$\Gamma_{1+2\gamma+l}\in H_{1+2\gamma+l}$\; such, that 
$$D\Gamma_{1+2\gamma+l}=f_{l}\;\;\text{ in }\;\R^{n+1}\;\;\text{ and } \;\;\lim_{y\rightarrow 0}
y^{1-2\gamma}\partial_{y}\Gamma_{1+2\gamma+l}=0\text{ on } \R^{n}\setminus 
\{0\}.$$
\end{lem}

\begin{pf}
This time we use
\begin{equation*}
\langle Q^{k}_{l}(y,x)=y^{2k}P_{l}(x)\mid \;\,k, \;l\in \N\;\;\text{ and }\;\;P_{l}\in 
\Pi_{l}\rangle \underset{\text{dense}}{\subset} C^{0}(\overline 
B_{1}^{\R^{n+1}}(0)\cap\R^{n+1}_{+}),
\end{equation*}
to obtain a orthogonal basis \;$E=\{e^{i}_{k}\}$\; for\; 
$L^{2}_{y^{1-2\gamma}}(S^{n}_{+})$
consisting of $D$-harmonics of the form 
\begin{equation*}
e^{i}_{k}=A_{m}\lfloor_{S^{n}_{+}}, \;\;\,A_{m}(y,x)=\sum 
y^{2l}P_{k-2l}(x), \;\;DA_{m}=0
\end{equation*}
and we have \;$
D_{S^{n}_{+}}e^{i}_{k}=k(k+n-2\gamma)y^{1-2\gamma}e^{i}_{k}.
$
Then for homogeneous \;$f, \;u$\; of degree \;$\lambda,\; \lambda+1+2\gamma$ solving
\begin{equation*}
\begin{cases}
Du=f\in L^{2}_{y^{2\gamma-1}}(\R^{n+1}_{+})\;\;\text{ in }\;\;\R^{n+1}_{+} \\
\lim_{y\rightarrow 0}y^{1-2\gamma}\partial_{y}u=0  \;\;\text{ on } \;\; \R^{n}
\end{cases}
\end{equation*}
is, when writing 
$
u=\sum a_{i,k}e^{i}_{k}\,,y^{2\gamma-1}f=\sum b_{j,l}e^{j}_{l},
$
equivalent to solving 
\begin{equation*}
\sum a_{i,k}(k(k+n-2\gamma)-(\lambda+1+2\gamma)(\lambda+n+1))e^{i}_{k}=\sum 
b_{j,l}e^{j}_{l}
\end{equation*}
and the latter system is always solvable in case
\begin{equation}\label{solvability_homogeneous_neumann}
k(k+n-2\gamma)-(\lambda+1+2\gamma)(\lambda+n+1)\neq 0 \,\,\text{ for all 
}\;\;k, \;n, \;\lambda \in \N. 
\end{equation}
As for proving the lemma there holds 
\begin{equation*}
deg(f_{l})=\lambda=m-n\;\;
\text{ and }\;\;
deg(e^{i}_{k})=k=m^{\prime}
\;\;\text{ for some }\,\;
m, \;m^{\prime}\in \N
\end{equation*}
and plugging this into 
\eqref{solvability_homogeneous_neumann} we verify for $\frac{1}{2}\neq \gamma \in (0,1)$
\begin{equation*}
m^{\prime}(m^{\prime}+n-2\gamma)-(m-n+1+2\gamma)(m+1)\neq 0\;\;\text{ for all } \;\;
n \;,m, \;m^{\prime}\in \N.
\end{equation*}
This shows homogeneous 
solvability, whereas regularity of the solution follows from Proposition 
\ref{green_regularity}.
\end{pf}
\vspace{6pt}

\noindent
Analogously to the case of the Poisson kernel  we may solve 
\eqref{green_defining_equation} successively using Lemma 
\ref{green_homogeneous_solvability_neumann}  and obtain
\begin{equation*}
\Gamma_{g}=\eta_{\xi}(\Gamma+\sum_{l=-n}^{m} \Gamma_{1+2\gamma+l})+\gamma_{m}
\end{equation*} 
for \;$m\geq 0$, where \;$\eta_{\xi}$\; is as in \eqref{etaxi} and a weak solution
\begin{equation*}
\begin{cases}
D_{g}\gamma_{m}=-D_{g}\left(\eta_{\xi}(\Gamma+\sum_{l=-n}^{m} 
\Gamma_{1+2\gamma+l})\right)=y^{1-2\gamma}h_{m}\;\; \text{ in }\;\;X \\
\lim_{y\rightarrow 0}y^{1-2\gamma}\partial_{y}\gamma_{m}=0 \,\;\text{ on } \;\;M 
\end{cases}
\end{equation*}
with \;$h_{m} \in C^{m,\alpha}$. As in the previous subsection a weak regularity statement is
sufficient for our purpose.
\begin{pro}
\label{green_regularity}$_{}$\\
Let \;$h\in y^{1-2\gamma}C^{2k+3,\alpha}(X)$\; and \;$u\in 
W^{1,2}_{y^{1-2\gamma}}(X)$\; be a weak solution of
\begin{equation*}
\begin{cases}
D_{g}u=h\;\;\text{ in }\;\;X \\
\lim_{y\rightarrow 0}y^{1-2\gamma}\partial_{y}u=0 \,\;\text{ on } \;\;M.
\end{cases}
\end{equation*}
Then \;$u$\; is of class \;$C^{2k,\beta}(X)$, provided \;$H_{g}=0$.
\end{pro}
\vspace{6pt}

\noindent
As in the previous subsection, putting these facts together before presenting the proof of Proposition \ref{green_regularity}, we  have the existence of \;$\Gamma_g$\; and can describe its asymptotics.
\begin{cor}
\label{Greens_function_asymptotics}$_{}$\\
Let \;$\frac{1}{2}\neq \gamma\in (0,1)$. Then \;$\Gamma_g$\; exists and we may expand in \;$g$-normal 
Fermi-coordinates around $\xi \in M$
\begin{equation*}
\Gamma_{g}(z,\xi)
\;\in
\eta_{\xi}(z)
\left( 
\frac{g_{n, \gamma}}{\vert z \vert^{n-2\gamma}}+\sum^{2m+3}_{l=-n}H_{1+2\gamma+l}(z)
\right)
+
C^{2m,\alpha}(X)
\end{equation*}
with \;$H_{l}\in C^{\infty}(\R^{n+1}_+\setminus \{0\})$\; being homogeneous of order 
\;$l$\; and \;$g_{n, \gamma}$\; is as in \eqref{gngamma}, provided \;$H_{g}=0$.
\end{cor}
\vspace{6pt}

\noindent
\begin{pfn}{ of Proposition \ref{green_regularity}}
\vspace{2pt}

\noindent
As in the previous subsection we use the Moser iteration argument. Indeed by exactly the same arguments as the ones used when proving Proposition 
\ref{poisson_regularity} we recover H\"older regularity \eqref{poisson_zamboni} 
and integrating the equation directly we find the analogue of
\eqref{poisson_equation_integrated_1}, namely
\begin{equation}\label{green_equation_integrated}
\begin{split}
\partial^{\alpha}_{x}u(y,x)
= &
u_{0}^{\alpha}(x) 
-
\int^{y}_{0}\sigma \tilde f^{\alpha}_{1}(\sigma,x)d\sigma
-
\int^{y}_{0}\sigma^{2\gamma-1}\int^{\sigma}_{0}\tau^{1-2\gamma} \tilde 
f^{\alpha}_{2}(\tau,x)d\tau d\sigma \\
= &
u^{\alpha}_{0}(x)+u^{\alpha}_{1}(y,x)+u^{\alpha}_{2}(y,x),
\end{split}
\end{equation}
where \;$\tilde f_{1}, \;\tilde f_{2}$\; are given by 
\eqref{poisson_tilde_f_structure}.
Let \;$h\in y^{1-2\gamma}C^{l,\lambda^{\prime}}$. Then 
\eqref{poisson_zamboni} and \eqref{poisson_tilde_f_structure} show
\begin{equation*}
\fa \vert \alpha\vert \leq l-2 \;:\,\;\tilde  
f_{i}^{\alpha}\in C^{0,\lambda}.
\end{equation*}
In particular \eqref{green_equation_integrated} implies
\begin{equation*}
\partial^{\alpha}_{x}u(y,x)
=
u_{0}^{\alpha}(x) +O(y),
\end{equation*} 
so \;$u_{0}^{\alpha}\in C^{l+2,\lambda}$\; anyway by interior regularity
and we may assume inductively
\begin{equation*}
\fa \vert \alpha\vert\leq  l-2-2m\; : \;\partial^{2m}_{y}\partial^{\alpha}_{x}u\in 
C^{0,\lambda}, 
\end{equation*}
whence according to \eqref{poisson_tilde_f_structure} 
\begin{equation*}
\fa \vert \alpha \vert\leq l-2-2(m+1)\; : \;\partial^{2m}_{y}\tilde 
f^{\alpha}_{i}\in C^{0,\lambda}. 
\end{equation*}
Then \eqref{green_equation_integrated} implies via Taylor expansion 
\begin{equation*}
\fa \vert \alpha \vert\leq 
l-2-2(m+1) \; : \;\partial_{y}^{2m+2}u_{i}^{\alpha}, \;\partial_{y}^{2m+2}\partial^{\alpha
}_{x}u \in C^{0,\lambda} 
\end{equation*}
Thus we have proven  \;$
\fa \vert \;\alpha \vert\leq 
l-2-2m \; : \;\partial_{y}^{2m}\partial^{\alpha
}_{x}u \in C^{0,\lambda} 
$
for some \;$\lambda>0$. However, since there are only even powers in 
the \;$y$-derivative, we only find \;$u\in C^{l-3,\lambda}$\; for \;$l\in 2\N$.
The proof is thereby complete.
\end{pfn}

\subsection{Green's function for the fractional conformal Laplacian }
In this short subsection we study the  Green's function \;$G_h^\gamma$\; of \;$P_h^\gamma$. We derive its existence and asymptotics as a consequence of the results of the previous subsections and formula \eqref{regreen}.

\begin{cor}
\label{Greens_function_asymptoticsgamma}$_{}$\\
Let \;$\frac{1}{2}\neq \gamma\in (0,1)$. Then \;$G_h$\; exists and we may expand in \;$h$-normal-coordinates
around $\xi \in M$
\begin{equation*}
G_{h}(x,\xi)
\;\in
\eta_{\xi}(x)
\left( 
\frac{g_{n, \gamma}}{\vert x \vert^{n-2\gamma}}+\sum^{2m+3}_{l=-n}H_{1+2\gamma+l}(x)
\right)
+
C^{2m,\alpha}(M)
\end{equation*}
with \;$H_{l}\in C^{\infty}(\R^{n}\setminus \{0\})$\; being homogeneous of order 
\;$l$, provided \;$H_{g}=0$.
\end{cor}
\vspace{10pt}

\noindent
To end this section, we give the proof of Theorem \ref{Greens_function_asymptoticsgamma1}.\\\\
\begin{pfn}{ of Theorem \ref{Greens_function_asymptoticsgamma1}}\\
It follows directly from Corollary \ref{Poisson_kernel_asymptotics}, Corollary \ref{Greens_function_asymptotics}, and Corollary \ref{Greens_function_asymptoticsgamma}.
\end{pfn}
\section{Locally flat conformal infinities of PE-manifolds}\label{pceinsten}
In this section we sharpen the results of Section \ref{fundso} in the case of Poincar\'e-Einstein manifold \;$(X, g^+)$\; with  locally flat conformal infinity \;$(M, [h])$. 
\subsection{Fermi-coordinates in this particular case}
 By our assumptions we have  
\begin{enumerate}[label=(\roman*)]
 \item a geodesic defining function \;$y$\; splitting the metric
\begin{equation*}
g=y^{2}g^{+}, \; g=dy^{2}+h_{y}\;\;\text{ near }\;\;M\;\;\text{ and }\;\; h=h_{y}\lfloor_{M}
\end{equation*} 
and for every \;$a\in M$\;  a conformal factor as in \eqref{conformal_factor_properties},
whose conformal metric \;$h_{a}=u_{a}^{\frac{4}{n-2\gamma}}h$\; close to \;$a$\; admits an Euclidean 
coordinate system, \;$h_{a}=\delta$\; on $B_{\epsilon}^{h_a}(a)$. As clarified in 
subsection \ref{confnonhom} and recalling Remark \ref{eq:minimal}, this gives rise to a geodesic defining 
function \;$y_{a}$, for which
\begin{equation}\label{dfmga}
g_{a}=y_{a}^{2}g^{+}, \;\; g_{a}=dy_{a}^{2}+h_{a,y_{a}}\;\;\text{near}\;\;M\;\;\text{ with }\;\; 
h_{a}=h_{a, y_a}\lfloor_{M}\;\;\text{ and }\;\;\delta=h_{a}\lfloor_{ 
B_{\epsilon}^{h_a}(a)},
\end{equation} 
the boundary \;$(M, [h_a])$\; is totally geodesic and the extension operator \;$D_{g_{a}}$\; is positive.

 \item
 as observed by Kim-Musso-Wei\cite{kmw1} in the case \;$n\geq 3$, cf. Lemma 43 in \cite{kmw1}, and for \;$n=2$\; due to Remark \ref{eq:minimal} and the existence of isothermal coordinates we have
\begin{equation}\label{flatness}
g_{a}=\delta+O(y_{a}^{n}) \;\;\text{ on } \;\;B_{\epsilon}^{g_a, +}(a)
\end{equation} 
in \;$g_{a}$-normal 
Fermi-coordinates around $a$for some small \;$\epsilon >0$. 
Therefore the previous results on the fundamental 
solutions in the case of an asymptotically hyperbolic manifold with minimal conformal infinity of Section \ref{fundso} are applicable. We collect them in the following subsection.
\end{enumerate}
\subsection{Fundamental solutions in this particular case}
In this subsection we sharpen the results of Section \ref{fundso} in the case of a Poincar\'e-Einstein manifold \;$(X, g^+)$\; with locally flat conformal infinity \;$(M, [h])$. 
\vspace{6pt}

\noindent
To do that  
let us first recall that  \;$K_{a}=K_{g_{a}}(\cdot, a)$, \;$\Gamma_{a}=\Gamma_{g_{a}}(\cdot, 
a)$ and \;$G_a=G_{h_a}(\cdot,a)$.
From \eqref{flatness} we then find  
\begin{equation*}
D_{g_{a}}K_{a}\in yH_{-2\gamma-2}C^{\infty},
\;\,\;\;
D_{g_{a}}\Gamma_{a}\in y^{1-2\gamma}H_{2\gamma-2}C^{\infty}
\end{equation*}
for the lowest order 
deficits in \eqref{poisson_local_equation_for_the_poisson_kernel} and 
\eqref{green_local_equation_for_the_greens_function}.
Then in view of Lemmas \ref{poisson_homogeneous_solvability_dirichlet}, \ref{green_homogeneous_solvability_neumann} the corresponding expansions given 
by Corollaries
\ref{Poisson_kernel_asymptotics}, \ref{Greens_function_asymptotics},  \ref{Greens_function_asymptoticsgamma}
are 
\begin{enumerate}[label=(\roman*)]
 \item \quad 
$
K_{a}(z)
\in
\eta_{\xi}(z)
\left(p_{n, \gamma} 
\frac{y^{2\gamma}}{\vert z 
\vert^{n+2\gamma}}+\sum^{2m+6}_{l=0}y^{2\gamma}H_{l-2\gamma}(z)
\right)
+
y^{2\gamma}C^{2m,\alpha}(X)
$ 

 \item \quad 
$
\Gamma_{a}(z)
\;\in
\eta_{\xi}(z)
\left( 
\frac{g_{n, \gamma}}{\vert z \vert^{n-2\gamma}}+\sum^{2m+4}_{l=0}H_{l+2\gamma}(z)
\right)
+
C^{2m,\alpha}(X)
$ 

 \item \quad 
$
G_{a}(x)
\;\in
\eta_{\xi}(x)
\left( 
\frac{g_{n, \gamma}}{\vert x \vert^{n-2\gamma}}+\sum^{2m+4}_{l=0}H_{l+2\gamma}(x)
\right)
+
C^{2m,\alpha}(M).
$ 

\end{enumerate}
Recalling \eqref{definition_H_l} there holds \;$y^{2\gamma}H_{l-2\gamma}\subset 
C^{m,\alpha}$\; for \;$l>m $\; and
$H_{l+2\gamma}\subset C^{m,\alpha}$\; for \;$l\geq m$. We have therefore proven the following result.
\begin{cor}
\label{cor_kernels_for_poincare_einstein_metrics}$_{}$\\
\noindent
Let \;$(X,  \;g^{+})$\; be a Poincar\'e-Einstein manifold with 
conformal infinity  \;$(M, [h])$\; of dimension \;$n= 2$\; or \;$n\geq 3$\; and \;$(M, [h])$\; is locally flat. If
\begin{equation*}
\frac{1}{2}\neq \gamma \in (0,1) \;\;\text{ and } \;\,\;
\lambda_{1}(-\Delta_{g^{+}})>s(n-s)\;\;
\text{ for } s=\frac{n}{2}+\gamma,
\end{equation*}
then the Poison kernel \;$K_{g}$\; and the Green's functions \;$\Gamma_{g}$\; and \;$G_{h}$\;  
respectively for
\begin{equation*}
\begin{cases}
D_{g}U=0 \;\;\text{ in }\;\; X \\
U=f \,\,\;\text{ on } \;\;M
\end{cases} 
\quad\quad
\begin{cases}
D_{g}U=0\;\; \text{ in } \;\;X \\
-d_{\gamma}^*\lim_{y\rightarrow 0}y^{1-2\gamma}\partial_{y}U=f \;\;\text{ on }\;\; M
\end{cases} 
\quad \text{and}\quad
\begin{cases}
P_h^\gamma u=f\;\;\text{ on }\;\;M
\end{cases}
\end{equation*}
are respectively of class \;$y^{2\gamma}C^{2,\alpha}$\; and \;$C^{2,\alpha}$\; away from the 
singularity and admit for every \;$a\in M$ locally in \;$g_{a}$-normal 
Fermi-coordinates an expansion around $a$
\begin{enumerate}[label=(\roman*)]
 \item \quad 
$
K_{a}(z)\;
\in
p_{n, \gamma}\frac{y^{2\gamma}}{\vert z \vert^{n+2\gamma}}
+
y^{2\gamma}H_{-2\gamma}(z)
+
y^{2\gamma}H_{1-2\gamma}(z)
+
y^{2\gamma}H_{2-2\gamma}(z)
+
y^{2\gamma}C^{2,\alpha}(X)
$

 \item \quad 
$
\Gamma_{a}(z)\;
\in
\frac{g_{n, \gamma}}{\vert z \vert^{n-2\gamma}}+H_{2\gamma}(z)+H_{1+2\gamma}(z)
+
C^{2,\alpha}(X)
$

 \item \quad 
$
G_{a}(x)\;
\in
\frac{g_{n, \gamma}}{\vert x \vert^{n-2\gamma}}+H_{2\gamma}(x)+H_{1+2\gamma}(x)
+
C^{2,\alpha}(M),
$

\end{enumerate}
where \;$g_a$\;is as in  \eqref{dfmga} and \;$H_{k}\in C^{\infty}(\overline{\R^{n}_{+}}\setminus \{0\})$\; are 
homogeneous of degree \;$k$. 
\end{cor}
\vspace{6pt}

\noindent
Finally, we give the proof of Theorem \ref{cor_kernels_for_poincare_einstein_metrics1}.\\\\
\begin{pfn}{ of Theorem \ref{cor_kernels_for_poincare_einstein_metrics1}}\\
It is exactly the statement of Corollary \ref{cor_kernels_for_poincare_einstein_metrics}.
\end{pfn}

\end{document}